\documentclass[12pt, draft]{article}
\usepackage{amssymb}
\usepackage{amsmath}
\usepackage{array}
\usepackage[T2A]{fontenc}
\usepackage[cp1251]{inputenc}
\newcommand{\Z}{\scriptstyle}
\newcommand{\D}{\displaystyle}
\newcommand{\prsum}{\mathop{{\sum}'}}
\newcommand{\stsum}{\mathop{{\sum}^{*}}}

\renewcommand{\S}{\mathhexbox278}
\renewcommand{\le}{\operatorname{\leqslant}}
\renewcommand{\ge}{\operatorname{\geqslant}}

\DeclareMathOperator{\RRe}{Re} \DeclareMathOperator{\IIm}{Im}
\DeclareMathOperator{\mmod}{mod}
\DeclareMathOperator{\vep}{\varepsilon}

\DeclareMathOperator{\vf}{\varphi}
\DeclareMathOperator{\vth}{\vartheta}

\multlinegap=0em \DeclareFontFamily{T1}{msb}{}
\DeclareFontShape{T1}{msb}{m}{ol}{<5> <6> <7> <8> <9> gen * msbm
<10> <10.95> <12> <14.4> <17.28> <20.74> <24.88> msbm10}{}
\DeclareSymbolFont{AMSb}{T1}{msb}{m}{ol} \multlinegap=0em
\begin{document}

\begin{center}
{\rmfamily\bfseries\normalsize On Karatsuba's Problem Concerning the
Divisor Function $\boldsymbol{\tau(n)}$\footnote{This research was
supported by the Programme of the President of the Russian
Federation `Young Candidates of the Russian Federation' (grant no.
МК-4052.2009.1).}}
\end{center}

\vspace{0.5cm}

\begin{center}
\textbf{M.A. Korolev}
\end{center}

\vspace{0.5cm}

\textbf{Abstract.} We study an asymptotic behavior of the sum
$\sum\limits_{n\le x}\frac{\D \tau(n)}{\D \tau(n+a)}$. Here
$\tau(n)$ denotes the number of divisors of $n$ and $a\ge 1$ is a
fixed integer.

\vspace{0.5cm}

\begin{flushleft}
\textbf{1. Introduction}
\end{flushleft}

\vspace{0.3cm}

This paper deals with a problem stated by A.A. Karatsuba in
November, 2004: to determine an asymptotic behavior of the following
sum
\[
S(x)\,=\,\sum\limits_{n\le x}\frac{\tau(n)}{\tau(n+1)}
\]
(here $\tau(n)$ denotes the number of positive divisors of $n$).
Since
\[
\frac{1}{x}\sum\limits_{n\le x}\tau(n)\,\sim\,\ln x,\quad
\frac{1}{x}\sum\limits_{n\le
x}\frac{1}{\tau(n)}\,\sim\,\frac{c}{\sqrt{\ln x}}
\]
for some $c>0$, the below assumption seems reasonable:
\[
\frac{1}{x}\,S(x)\,\sim\,\frac{K}{\sqrt{\ln x}}\cdot \ln
x\,\sim\,K\sqrt{\ln x},\quad K>0.
\]
The aim of this paper is to prove the following

\vspace{0.2cm}

\textbf{Theorem.} \emph{Let} $a$ \emph{be a fixed integer}, $a\ge
1$, \emph{and let}
\[
S_{a}(x)\,=\,\sum\limits_{n\le x}\frac{\tau(n)}{\tau(n+a)}.
\]
\emph{Then}
\[
S_{a}(x)\,=\,K(a)x\sqrt{\ln x}\,+\,O(x\ln\ln x),
\]
\emph{as} $x\to + \infty$. \emph{The constant} $K(a)$ \emph{has the
form} $K(a) = K\cdot\kappa(a)$, \emph{where}
\begin{align*}
& K =
\frac{1}{\sqrt{\pi}}\prod\limits_{p}\Bigl(\frac{1}{\sqrt{p(p-1)}}\,+\,\sqrt{1
-
\frac{1}{p}}(p-1)\ln\frac{p}{p-1}\Bigr)\,=\,0.757\,827\,651\ldots,\\
& \kappa(a) = \beta(a)\prod\limits_{p|a}\frac{\D 1 + \sum\limits_{k
= 1}^{+\infty}\frac{\D e_{a}(p^{k})}{\D k+1}\,p^{-k}}{\D 1 +
\beta(p)\sum\limits_{k = 1}^{+\infty}\frac{p^{-k}}{k+1}},\\
& e_{a}(n) =
\frac{1}{\beta(a)}\sum\limits_{d|(a,n)}\beta\Bigl(\frac{an}{d^{2}}\Bigr),\quad
\beta(a) = \prod\limits_{p|a}\frac{(p-1)^{2}}{p^{2}-p+1}.
\end{align*}

\vspace{0.2cm}

For the below, we need the following notations:

\vspace{0.2cm}

$\vf(q)$ denotes Euler function;

$\chi$ denotes Dirichlet's character modulo $q$, $q\ge 3$;

the symbols $\sum\limits_{\substack{\chi\,\text{mod}\, q \\
\chi\ne \chi_{0}}}$ and $\sum\limits_{\chi\ne \chi_{0}}$ denote the
sums over all non-principal characters modulo $q$;

the symbols $\stsum\limits_{\chi\,\text{mod}\, q}$ and
$\stsum\limits_{\chi}$ denote the sums over all primitive characters
modulo $q$;

$s = \sigma + it$, $\sigma, t$ are real numbers;

$L(s,\chi)$ is Dirichlet's $L$-function corresponding to the
character $\chi$;

the symbol $N(\sigma; T,\chi)$ means the number of zeros of
$L(s,\chi)$ in the rectangle $\sigma < \RRe s \le 1$, $|\IIm s|\le
T$;

$(a,b)$ denotes the great common divisor of $a$ and $b$;

symbols $\theta, \theta_{1}, \theta_{2}, \ldots$ denote complex
numbers such that $|\theta|$,$|\theta_{1}|$,$|\theta_{2}|$,$\ldots
\le 1$, in general, different in different relations.

\vspace{0.2cm}

\begin{flushleft}
\textbf{2. Auxiliary assertions}
\end{flushleft}

\vspace{0.2cm}

\textbf{Lemma 1.} \emph{Suppose} $S(t)$ \emph{is a smooth
complex-valued function for} $t_{0}\le t\le t_{k}$,
$t_{0}<t_{1}<\ldots < t_{k}$ \emph{and} $\D \min_{0\le j\le
k-1}(t_{j+1} - t_{j}) = \delta > 0$. \emph{Then the following
inequality holds}:
\[
\sum\limits_{j =
1}^{k}|S(t_{j})|^{2}\,\le\,\frac{1}{\delta}\,j_{1}\,+\,2\sqrt{j_{1}\,j_{2}},
\]
\emph{where}
\[
j_{1}\,=\,\int_{t_{0}}^{t_{k}}|S(t)|^{2}dt,\quad
j_{2}\,=\,\int_{t_{0}}^{t_{k}}|S'(t)|^{2}dt.
\]
For the proof, see \cite[Chapter VII, \S 1]{Karatsuba_1983}.

\vspace{0.2cm}

\textbf{Lemma 2.} \emph{Suppose} $M, N, Q$ \emph{are integers}.
\emph{Then for any sequence of complex numbers} $a_{n}$ \emph{the
following estimation holds}:
\[
\sum\limits_{q\le Q}\stsum\limits_{\chi\,\text{mod}\,
q}\Bigl|\sum\limits_{n =
M+1}^{M+N}a_{n}\chi(n)\Bigr|^{2}\,\ll\,(Q^{2}\,+\,N)\sum\limits_{n =
M+1}^{M+N}|a_{n}|^{2},
\]
\emph{where the constant in the symbol} $\ll$ \emph{is absolute}.

For the proof, see \cite[Chapter IX, answers to
problems]{Karatsuba_1983}.

\vspace{0.2cm}

\textbf{Lemma 3.} \emph{Suppose} $q\ge 3$, $\chi$ \emph{is
non\,-principle Dirichlet's character modulo} $q$. \emph{Then for
any} $s, Y$ \emph{such that} $Re s \ge \sigma_{0}>0$ \emph{and}
$Y\ge q(|t|+1)/\pi$ \emph{the following equality holds}:
\[
L(s,\chi)\,=\,\sum\limits_{n\le
Y}\frac{\chi(n)}{n^{s}}\,+\,O\bigl(qY^{-\sigma}\bigr),
\]
\emph{where the constant in the symbol} $O$ \emph{is absolute}.

\vspace{0.2cm}

For the proof, see \cite[\S 26]{Davenport_1971}.

\vspace{0.2cm}

\textbf{Lemma 4.} \emph{Let} $c$ \emph{be a sufficiently small
positive absolute constant and let} $q\ge 3$. \emph{Suppose} $\chi$
\emph{is a complex character modulo} $q$. \emph{Then the function}
$L(s,\chi)$ \emph{has no zeros in the domain}
\[
\RRe s \,>\,1\,-\,\frac{c}{\ln{q(|t|+1)}},\quad -\infty < t <
+\infty;
\]
\emph{now if} $\chi$ \emph{is a real non-principal Diriclet's
character modulo} $q$ \emph{then the function} $L(s,\chi)$ \emph{has
no zeros in the domain}
\[
\RRe s \,>\,1\,-\,\frac{c}{\ln{q(|t|+1)}},\quad |t|>0.
\]
\vspace{0.2cm}

For the proof, see \cite[Chapter IX, \S 2]{Karatsuba_1983}.

\vspace{0.2cm}

\textbf{Lemma 5} (Siegel). \emph{For any} $\vep$,
$0<\vep<\tfrac{1}{2}$, \emph{there exists} $c = c(\vep) > 0$
\emph{such that if} $\chi$ \emph{is a real character modulo} $q$
\emph{and} $\beta$ \emph{is a real zero of} $L(s,\chi)$, \emph{then}
\[
\beta\,<\,1\,-\,\frac{c}{q^{\vep}}.
\]
\vspace{0.2cm}

For the proof, see \cite[Chapter IX, \S 2]{Karatsuba_1983}. The
constant $c = c(\vep)$ is not effective. This means that it is
impossible to find or estimate $c(\vep)$ from a given $\vep$.
Therefore all statements (including main theorem of this paper) in
which this lemma is essentially used are ineffective, too.

\vspace{0.2cm}

\textbf{Lemma 6} (Montgomery). \emph{For any} $Q\ge 3$, $T\ge 3$
\emph{the following estimation holds}:
\[
\sum\limits_{q\le Q}\stsum\limits_{\chi\,\text{mod}\,q}N(\sigma;
T,\chi)\,\ll\,\bigl(Q^{2}T\bigr)^{\vth(\sigma)}(\ln QT)^{14}
\]
\emph{where}
\begin{equation*}
\vth(\sigma)\,=\,
\begin{cases}
\frac{\D 3(1-\sigma)}{\D 2-\sigma\mathstrut}, & \text{if}\quad
\tfrac{1}{2}\le
\sigma \le \tfrac{4}{5},\\
\frac{\D 2(1-\sigma)\mathstrut}{\D \sigma\mathstrut}, &
\text{if}\quad \tfrac{4}{5}\le \sigma \le 1,
\end{cases}
\end{equation*}
\emph{and the constant in the symbol} $\ll$ \emph{is absolute}.
\vspace{0.2cm}

For the proof, see \cite[Chapter 12]{Montgomery_1974}.

\vspace{0.2cm}

\textbf{Lemma 7.} \emph{Suppose} $q\ge 1$ \emph{is an integer},
$\chi$ \emph{is a character modulo} $q$ \emph{and} $\varrho = \beta
+ i\gamma$ \emph{runs through all non-trivial zeros of} $L(s,\chi)$.
\emph{Then}
\[
\sum\limits_{|\gamma|\le T}\frac{1}{|\gamma|+1}\,\ll\,(\ln qT)^{2},
\]
\emph{as} $T\to +\infty$.

\vspace{0.2cm}

This estimation follows from asymptotic formula for $N(T,\chi)$ -
the number of zeros of $L(s,\chi)$ in the rectangle $0\le \RRe  s\le
1$, $|\IIm s|\le T$.

\pagebreak

\begin{flushleft}
\textbf{3. Basic assertions}
\end{flushleft}

\vspace{0.2cm}

\textbf{Lemma 8.} \emph{Let} $\chi$ \emph{be a non-primitive
character modulo} $q = q_{1}r$ \emph{induced by primitive character}
$\chi_{1}$ \emph{modulo} $q_{1}$. \emph{Then}
\[
|L(s,\chi)|\,\le\,\tau(r)|L(s,\chi_{1})|
\]
\emph{in the half-plane} $\RRe s\ge 0$.

\vspace{0.2cm}

\emph{Proof.} Using the formula
\[
L(s,\chi)\,=\,L(s,\chi_{1})\prod\limits_{p|q,\,p\nmid q_{1}}\Bigl(1
- \frac{\chi_{1}(p)}{p^{s}}\Bigr),
\]
we get
\[
|L(s,\chi)|\,\le\,|L(s,\chi_{1})|\,\prod\limits_{p|r}\Bigl(1 +
\frac{1}{p^{\sigma}}\Bigr)\,\le\,|L(s,\chi_{1})|
\prod\limits_{p|r}2\,\le\,\tau(r)|L(s,\chi_{1})|.
\]

\vspace{0.2cm}

\textbf{Lemma 9.} \emph{Suppose} $f(n)$ \emph{is non-negative
multiplicative function such that} $f(n) = O(n^{\vep})$, $0<\vep <
\tfrac{1}{2}$, \emph{and the function}
\[
F(s) = \sum\limits_{n = 1}^{+\infty}\frac{f(n)}{n^{s}}
\]
\emph{satisfies the identity}
\[
F(s)\,=\,\sqrt{\zeta(s)}\,\Phi(s)
\]
\emph{for} $\RRe s > 1$ ($\sqrt{z} > 0$ \emph{for} $z > 0$),
\emph{where} $\Phi(s)$ \emph{is regular in the half-plane} $\RRe s >
\tfrac{1}{2}$ \emph{and obeys the estimate}
\[
|\Phi(\sigma + it)|\,\ll\, \max{\Bigl\{1,\bigl(\sigma -
\tfrac{1}{2}\bigr)^{-c}\Bigr\}}
\]
\emph{for any} $\sigma > \tfrac{1}{2}$ \emph{and some} $c > 0$.
\emph{Then}
\[
\sum\limits_{n\le x}f(n)\,=\,\frac{x}{\sqrt{\ln
x}}\biggl(\frac{\Phi(1)}{\sqrt{\pi}}\,+\,O\Bigl(\frac{1}{\ln
x}\Bigr)\biggr)
\]
\emph{where the constant in the symbol} $O$ \emph{depends on} $f$
\emph{only}.

\vspace{0.2cm}

\emph{Proof.} Since $f(n) = O(n^{\vep})$, without loss of generality
we may assume that $x$ has the form $N + \tfrac{1}{2}$ for some
integer $N$. Suppose $T$ differs from the imaginary part of any zero
of $\zeta(s)$ and $2\le T \le x$. Then, by Perron's formula we get
\[
\sum\limits_{n\le x}f(n)\,=\,I\,+\,O\Bigl(\frac{x}{T}\ln x\Bigr)
\]
where
\[
I\,=\,\frac{1}{2\pi
i}\int_{b-iT}^{b+iT}F(s)\,\frac{x^{s}}{s}\,ds,\quad b = 1 +
\frac{1}{\ln x}.
\]
Let $c_{1}>0$ be small positive constant such that $\zeta(s)$ has no
zeros in the rectangle with vertices $1\pm iT$, $\alpha \pm iT$,
\[
\alpha\,=\, 1\,-\,c_{1}(\ln T)^{-\frac{\Z 2}{\Z 3\mathstrut}}(\ln\ln
T)^{-\frac{\Z 1}{\Z 3\mathstrut}}.
\]
Then, by the identity $F(s)=\sqrt{\zeta(s)}\,\Phi(s)$ the function
$F(s)$ continues analytically to the domain $\alpha\le \RRe s \le
1$, $|\IIm s|\le T$ with horizontal cut going straight from the
point $s = \alpha$ to the point $s = 1$. By Cauchy's theorem,
\[
I\,=\,-\sum\limits_{k = 1}^{6}I_{k},
\]
where the symbols $I_{1},\ldots, I_{4}$ denote the integrals along
the segments connecting points $b+iT$, $\alpha + iT$, $\alpha$,
$\alpha - iT$, $b - iT$ and the symbols $I_{5}, I_{6}$ denote the
integrals along the upper and lower edges of the cut respectively.

Since the bound
\[
\zeta(\sigma + it)\,=\,O\bigl(\ln^{\frac{\Z 2}{\Z
3\mathstrut}}(|t|+2)\bigr)
\]
holds along the contour, we obtain:
\begin{align*}
& |I_{1}| + |I_{4}|\,\ll\, \int_{\alpha}^{b}\bigl(|F(\sigma +
iT)|\,+\,|F(\sigma -
iT)|\bigr)\frac{x^{\sigma}d\sigma}{\sqrt{\sigma^{2}+T^{2}}}\,\ll\,\frac{x}{T}(\ln
T)^{\frac{\Z 1}{\Z 3\mathstrut}},\\
& |I_{2}| + |I_{3}| \ll \int_{-T}^{T}|F(\alpha +
it)|\frac{x^{\alpha}\,dt}{\sqrt{\alpha^{2}+t^{2}}}\ll
x^{\alpha}\int_{0}^{T}\frac{\ln^{\frac{\Z 1}{\Z
3\mathstrut}}(t+2)}{t+2}dt\ll x^{\alpha}(\ln x)^{\frac{\Z 4}{\Z
3\mathstrut}}.
\end{align*}
Consider the function $u(s) = (s-1)\zeta(s)$, $u(1) = 1$. Since
$u(s)\ne 0$ for $|s-1|\le \tfrac{1}{2}$, it follows that
\begin{align*}
& I_{5}+I_{6} = \frac{1}{2\pi
i}\int_{\alpha}^{1}\bigl(\sqrt{\zeta(\sigma + i\cdot
0)}\,-\,\sqrt{\zeta(\sigma - i\cdot
0)}\bigr)\Phi(\sigma)\,\frac{x^{\sigma}}{\sigma}d\sigma\,=\\
& =\,\frac{1}{2\pi
i}\int_{\alpha}^{1}\Bigl(\frac{1}{\sqrt{\sigma-1+i\cdot 0}} -
\frac{1}{\sqrt{\sigma-1-i\cdot
0}}\Bigr)\Phi(\sigma)\sqrt{u(\sigma)}\,\frac{x^{\sigma}}{\sigma}d\sigma\,=\\
& =\,\frac{1}{2\pi
i}\int_{0}^{1-\alpha}\Bigl(\frac{1}{i\sqrt{v}}\,-\,\frac{1}{(-i\sqrt{v})}\Bigr)\,\frac{\Phi(1-v)\sqrt{u(1-v)}}{1-v}\,x^{1-v}\,dv\,=\\
&
=\,-\,\frac{x}{\pi}\int_{0}^{1-\alpha}\frac{x^{-v}}{\sqrt{v}}\,\frac{\Phi(1-v)\sqrt{u(1-v)}}{1-v}\,dv.
\end{align*}
The functions $\Phi(s)$ and $u(s)$ are bounded for
$|s-1|\le\tfrac{1}{3}$. This implies that
\[
\frac{\Phi(1-v)\sqrt{u(1-v)}}{1-v}\,=\,\Phi(1)\,+\,O(|v|)
\]
where the constant in the $O$\,-symbol depends on $\Phi$ only. Thus
we obtain
\begin{align*}
& I_{5} + I_{6} =
-\,\frac{x}{\pi}\int_{0}^{1-\alpha}\frac{x^{-v}}{\sqrt{v}}\bigl(\Phi(1)\,+\,O(v)\bigr)dv\,=\\
&
=\,-\,\frac{x}{\pi}\,\Phi(1)\Bigl(\int_{0}^{+\infty}\frac{x^{-v}}{\sqrt{v}}\,dv\,-\,\int_{1-\alpha}^{+\infty}\frac{x^{-v}}{\sqrt{v}}\,dv\Bigr)
\,+\,O\Bigl(x\int_{0}^{+\infty}\sqrt{v}x^{-v}\,dv\Bigr).
\end{align*}
Since
\[
(1-\alpha)\ln x = c_{1}(\ln x)(\ln T)^{-\frac{\Z 2}{\Z
3\mathstrut}}(\ln\ln T)^{-\frac{\Z 1}{\Z 3\mathstrut}}\ge
c_{1}\Bigl(\frac{\ln x}{\ln\ln x}\Bigr)^{\frac{\Z 1}{\Z
3\mathstrut}} > 1,
\]
we get
\begin{align*}
&
\int_{1-\alpha}^{+\infty}\frac{x^{-v}}{\sqrt{v}}\,dv\,=\,\frac{1}{\sqrt{\ln
x}}\int_{(1-\alpha)\ln x}^{+\infty}\frac{e^{-w}}{\sqrt{w}}\,dw <
\frac{1}{\sqrt{\ln x}}\int_{(1-\alpha)\ln x}^{+\infty}e^{-w}\,dw =
\frac{x^{\alpha - 1}}{\sqrt{\ln x}},\\
& I_{5}+I_{6} = -\,\frac{\Phi(1)}{\sqrt{\pi}}\,\frac{x}{\sqrt{\ln
x}}\,+\,O\Bigr(\frac{x^{\alpha - 1}}{\sqrt{\ln x}}\Bigr) +
O\Bigl(\frac{x}{(\ln x)^{3/2}}\Bigr).
\end{align*}
Using the above inequalities for $I_{1},\ldots, I_{4}$ and
substituting $T = e^{\sqrt{\ln x}}$, we finally obtain
\[
I\,=\,\frac{\Phi(1)}{\sqrt{\pi}}\,\frac{x}{\sqrt{\ln
x}}\Bigl(1\,+\,O\Bigl(\frac{1}{\ln x}\Bigr)\Bigr).
\]
The proof is complete.

\vspace{0.2cm}

\textbf{Lemma 10.} \emph{Let} $m$ \emph{be an integer}, $m\ge 1$.
\emph{Then}
\[
\sum\limits_{\substack{q\le x \\ (q,m) =
1}}\frac{1}{\vf(q)}\,=\,C\beta(m)\Bigl(\ln
x\,+\,\gamma\,-\,\sum\limits_{p}\frac{\ln p}{p^{2}-p+1} +
\sum\limits_{p|m}\frac{p^{2}\ln p}{(p-1)(p^{2}-p+1)}\Bigr)\,+
\]
\[
+\,O\Bigl(\frac{\ln^{2}x}{x}\Bigr)\,+\,O\Bigl(\frac{\tau(m)\ln
x}{x}\Bigr),
\]
\emph{where} $\gamma$ \emph{is Euler's constant},
\[
C\,=\,\prod\limits_{p}\Bigl(1\,+\,\frac{1}{p(p-1)}\Bigr)\,=\,\frac{\zeta(2)\zeta(3)}{\zeta(6)},
\]
\emph{and the constants in} $O$'\emph{s are absolute}.

\vspace{0.2cm}

\emph{Proof.} Note that
\[
\frac{1}{\vf(q)} =
\frac{1}{q}\prod\limits_{p|q}\Bigl(1\,-\,\frac{1}{p}\Bigr)^{-1}\!\!
= \frac{1}{q}\prod\limits_{p|q}\Bigl(1\,+\,\frac{1}{p-1}\Bigr) =
\frac{1}{q}\prod\limits_{p|q}\Bigl(1\,+\,\frac{1}{\vf(p)}\Bigr) =
\frac{1}{q}\sum\limits_{d|q}\frac{\mu^{2}(d)}{\vf(d)}.
\]
Let us use the prime sign for the summation over the numbers coprime
to $m$. Thus we get
\begin{align*}
&\prsum\limits_{q\le x}\frac{1}{\vf(q)} = \prsum\limits_{q\le
x}\frac{1}{q}\sum\limits_{d|q}\frac{\mu^{2}(d)}{\vf(d)} =
\prsum\limits_{d\le
x}\frac{\mu^{2}(d)}{\vf(d)}\prsum\limits_{\substack{q\le x\\ q\equiv
0 (\mmod d)}}\frac{1}{q}\,=\\
& =\,\prsum\limits_{d\le
x}\frac{\mu^{2}(d)}{\vf(d)}\prsum\limits_{k\le
x/d}\frac{1}{kd}\,=\,\prsum\limits_{d\le
x}\frac{\mu^{2}(d)}{d\vf(d)}\prsum\limits_{k\le
x/d}\frac{1}{k}\,=\\
& =\,\prsum\limits_{d\le
x}\frac{\mu^{2}(d)}{d\vf(d)}\sum\limits_{k\le
x/d}\Bigl(\sum\limits_{\delta |
(k,m)}\mu(\delta)\Bigr)\frac{1}{k}\,=\,\prsum\limits_{d\le
x}\frac{\mu^{2}(d)}{d\vf(d)}\sum\limits_{\delta|
m}\mu(\delta)\sum\limits_{\substack{k\le x/d \\ k\equiv 0 (\mmod
\delta)}}\frac{1}{k}\,=\\
& =\,\sum\limits_{\delta| m}\mu(\delta)\prsum\limits_{d\le
x}\frac{\mu^{2}(d)}{d\vf(d)}\sum\limits_{r\le
x/(d\delta)}\frac{1}{r\delta}\,=\,\\
& =\,\sum\limits_{\delta | m}
\frac{\mu(\delta)}{\delta}\prsum\limits_{d\le
x}\frac{\mu^{2}(d)}{d\vf(d)}\Bigl(\ln\frac{x}{d\delta}\,+\,\gamma\,+\,O\Bigl(\frac{d\delta}{x}\Bigr)\Bigr).
\end{align*}
The contribution of the error term if the brackets to the initial
sum does not exceed in order
\[
\sum\limits_{\delta | m}\frac{1}{\delta}\prsum\limits_{d\le
x}\frac{1}{d\vf(d)}\,\frac{d\delta}{x}\,=\,\frac{1}{x}\Bigl(\sum\limits_{\delta
| m}1\Bigr)\sum\limits_{d\le
x}\frac{1}{\vf(d)}\,=\,O\Bigl(\frac{\tau(m)\ln x}{x}\Bigr).
\]
The contribution of all other terms has the form
\begin{align*}
& \sum\limits_{\delta |
m}\frac{\mu(\delta)}{\delta}\prsum\limits_{d\le
x}\frac{\mu^{2}(d)}{d\vf(d)}\,(\ln x + \gamma - \ln d - \ln
\delta)\,=\,\\
& =\, (\ln x + \gamma)\Bigl(\sum\limits_{\delta |
m}\frac{\mu(\delta)}{\delta}\Bigr)\prsum\limits_{d\le
x}\frac{\mu^{2}(d)}{d\vf(d)}\,-\,\Bigl(\sum\limits_{\delta |
m}\frac{\mu(\delta)}{\delta}\Bigr)\prsum\limits_{d\le
x}\frac{\mu^{2}(d)\ln d}{d\vf(d)}\;-\\
& -\, \Bigl(\sum\limits_{\delta |
m}\frac{\mu(\delta)\ln\delta}{\delta}\Bigr)\prsum\limits_{d\le
x}\frac{\mu^{2}(d)}{d\vf(d)}.
\end{align*}
Let us replace all the sums over $d\le x$ by infinite sums over $d$
coprime to $m$. Using the inequality
\[
\sum\limits_{q\le x}\frac{1}{\vf(q)}\,\ll\,\ln x,
\]
we obtain:
\begin{align*}
& \prsum\limits_{d\le
x}\frac{\mu^{2}(d)}{d\vf(d)}\,=\,\prsum\limits_{d =
1}^{+\infty}\frac{\mu^{2}(d)}{d\vf(d)}\,+\,O\Bigl(\frac{\ln
x}{x}\Bigr),\\
& \prsum\limits_{d\le x}\frac{\mu^{2}(d)\ln
d}{d\vf(d)}\,=\,\prsum\limits_{d=1}^{+\infty}\frac{\mu^{2}(d)\ln
d}{d\vf(d)}\,+\,O\Bigl(\frac{\ln^{2}x}{x}\Bigr).
\end{align*}
Since
\begin{align*}
& \prsum\limits_{d=1}^{+\infty}\frac{\mu^{2}(d)\ln
d}{d\vf(d)}\,=\,\prod\limits_{p\nmid
m}\Bigl(1\,+\,\frac{1}{p(p-1)}\Bigr)\,=\,\prod\limits_{p}\Bigl(1\,+\,\frac{1}{p(p-1)}\Bigr)\prod\limits_{p
| m}\frac{p(p-1)}{p^{2}-p+1}\,=\\
& =\, C\prod\limits_{p | m}\frac{p(p-1)}{p^{2}-p+1},
\end{align*}
we have
\begin{align*}
& \prsum\limits_{n\le x}\frac{1}{\vf(q)} = C\,\frac{\vf(m)}{m}\,(\ln
x + \gamma)\prod\limits_{p |
m}\frac{p(p-1)}{p^{2}-p+1}\,-\,\frac{\vf(m)}{m}\prsum\limits_{d=1}^{+\infty}\frac{\mu^{2}(d)\ln
d}{d\vf(d)}\,-\\
& -\, C\prod\limits_{p |
m}\frac{p(p-1)}{p^{2}-p+1}\Bigl(\sum\limits_{\delta |
m}\frac{\mu(\delta)\ln \delta}{\delta}\Bigr).
\end{align*}
Note that non-zero terms in the sum over $\delta| m$ correspond to
squarefree divisors $\delta$, so we have
\[
\ln\delta\,=\,\sum\limits_{p | \delta}\ln p.
\]
Therefore,
\begin{align*}
& \sum\limits_{\delta | m}\frac{\mu(\delta)\ln
\delta}{\delta}\,=\,\sum\limits_{\delta |
m}\frac{\mu(\delta)}{\delta}\sum\limits_{p | \delta}\ln
p\,=\,\sum\limits_{p|m}\ln p \sum\limits_{\substack{\delta | m \\
\delta\equiv 0 (\mmod p)}}\frac{\mu(\delta)}{\delta}\,=\\
& =\,\sum\limits_{p | \delta}\ln p\sum\limits_{\delta_{1}| \frac{\Z
m}{\Z p},\;(\delta_{1},p) =
1}\frac{\mu(p)\mu(\delta_{1})}{p\delta_{1}}\,=\,-\sum\limits_{p|m}\frac{\ln
p}{p}\sum\limits_{\delta_{1}|\frac{\Z m}{\Z p},\;(\delta_{1},p) =
1}\frac{\mu(\delta_{1})}{\delta_{1}}\,=\\
& =\,-\,\sum\limits_{p|m}\frac{\ln p}{p}\prod\limits_{r\mid
\frac{m}{p},\;r\ne
p}\Bigl(1\,-\,\frac{1}{r}\Bigr)\,=\,-\,\sum\limits_{p|m}\frac{\ln
p}{p}\cdot\frac{\D \prod\limits_{p|m}\Bigl(1 - \frac{1}{r}\Bigr)}{\D
1 - \frac{1}{p}\mathstrut}\,=\\
& =\,-\,\frac{\vf(m)}{m}\sum\limits_{p|m}\frac{\ln p}{p-1}.
\end{align*}
Next, by the same arguments we obtain
\begin{align*}
& \prsum\limits_{d = 1}^{+\infty}\frac{\mu^{2}(d)\ln
d}{d\vf(d)}\,=\,\prsum\limits_{d =
1}^{+\infty}\frac{\mu^{2}(d)}{d\vf(d)}\sum\limits_{p|d}\ln
p\,=\,\prsum\limits_{p}\ln p\prsum\limits_{\substack{d = 1
\\ d\equiv 0 (\mmod p)}}^{+\infty}\frac{\mu^{2}(d)}{d\vf(d)}\,=\\
& =\,\prsum\limits_{p}\ln p\prsum\limits_{\substack{d_{1} = 1
\\ (d_{1},p) =
1}}^{+\infty}\frac{\mu^{2}(pd_{1})}{pd_{1}\vf(pd_{1})}\,=\,\prsum\limits_{p}\frac{\ln
p}{p\vf(p)}\sum\limits_{\substack{d_{1} = 1
\\ (d_{1},\,pm) =
1}}^{+\infty}\frac{\mu^{2}(d_{1})}{d_{1}\vf(d_{1})}\,=\\
& =\,\prsum\limits_{p}\frac{\ln p}{p(p-1)}\prod\limits_{q\nmid
pm}\Bigl(1\,+\,\frac{1}{q(q-1)}\Bigr)\,=\\
& =\,\prsum\limits_{p}\frac{\ln p}{p(p-1)}\cdot\frac{\D
\prod\limits_{q}\Bigl(1 + \frac{1}{q(q-1)}\Bigr)}{\D \Bigl(1 +
\frac{1}{p(p-1)}\Bigr)\prod\limits_{q|m}\Bigl(1 +
\frac{1}{q(q-1)}\Bigr)}\,=\\
& =\,C\prsum\limits_{p}\frac{\ln
p}{p(p-1)}\cdot\frac{p(p-1)}{p^{2}-p+1}\prod\limits_{q|m}\Bigl(1+\frac{1}{q(q-1)}\Bigr)^{\!-1}\,=\\
& =\,C\prsum\limits_{p}\frac{\ln
p}{p^{2}-p+1}\prod\limits_{q|m}\frac{q(q-1)}{q^{2}-q+1}\,=\\
&
=\,C\prod\limits_{q|m}\frac{q(q-1)}{q^{2}-q+1}\,\biggl(\sum\limits_{p}\frac{\ln
p}{p^{2}-p+1}\,-\,\sum\limits_{p|m}\frac{\ln p}{p^{2}-p+1}\biggr).
\end{align*}
Substituting these relations in the above formula for the initial
sum and taking into account the relation
\[
\frac{\vf(m)}{m}\prod\limits_{p|m}\frac{p(p-1)}{p^{2}-p+1}\,=\,\prod\limits_{p|m}\frac{(p-1)^{2}}{p^{2}-p+1}\,=\,\beta(m),
\]
we get
\begin{align*}
& \prsum\limits_{q\le x}\frac{1}{\vf(q)}\,=\\
& =\,C\beta(m)\biggl(\ln x + \gamma - \sum\limits_{p}\frac{\ln
p}{p^{2}-p+1} + \sum\limits_{p|a}\frac{\ln p}{p^{2}-p+1} +
\sum\limits_{p|a}\frac{\ln p}{p-1}\biggr)\,+\\
& +\,O\Bigl(\frac{\ln^{2}x}{x}\Bigr)\,+\,O\Bigl(\frac{\tau(a)\ln
x}{x}\Bigr).
\end{align*}

\vspace{0.2cm}

\textbf{Corollary.} \emph{Under the same conditions},
\[
\prsum\limits_{q\le x}\frac{1}{\vf(q)}\,=\,C\beta(m)\ln
x\,+\,O\bigl(\tau(m)\bigr)
\]
\emph{where the constant in the symbol} $O$ \emph{is absolute}.

\vspace{0.2cm}

\textbf{Lemma 11.} \emph{Let} $d\ge 2$ \emph{be a fixed integer.
Suppose} $\delta$ \emph{runs over an increasing sequence that
contains} $1$ \emph{and all integers which are not divisible by
prime numbers coprime to} $d$. \emph{Suppose also}
\[
D_{1}(x)\,=\,\sum\limits_{\delta\le x}1,\quad
D_{2}(x)\,=\,\sum\limits_{\delta > x}\frac{1}{\delta}.
\]
\emph{Then}
\[
D_{1}(x)\,\ll\,(\ln x)^{s}\quad D_{2}(s)\,\ll\,\frac{(\ln
x)^{s}}{x},
\]
\emph{where} $1\le s\le\tau(d)$.

\vspace{0.2cm}

\emph{Proof.} Let $p_{1}^{\alpha_{1}}\ldots p_{s}^{\alpha_{s}}$ be
the unique decomposition of $d$ into prime-powers. Then $D_{1}(x)$
equals to the number of solutions of the inequality
$p^{\beta_{1}}\ldots p_{s}^{\beta_{s}}\le x$ or
\[
\beta_{1}\ln p_{1}\,+\,\ldots\,+\,\beta_{s}\ln p_{s}\,\le\,\ln x
\]
with non-negative integers $\beta_{1},\ldots, \beta_{s}$. Since $\ln
p_{1}\ge \ln 2$, $\ldots$, $\ln p_{s}\ge \ln 2$, it follows that
$D_{1}(x)$ does not exceed the number of solutions of the inequality
\[
\beta_{1}\,+\,\ldots\,+\,\beta_{s}\,\le m,\quad m\,=\,
\left[\frac{\ln x}{\ln 2}\right],
\]
that is $D_{1}(x)\le \binom{\D m+s}{\D s}$. Applying Stirling's
formula and Cauchy's inequality we obtain
\begin{align*}
& \binom{\D m+s}{\D s}\,=\,\frac{1}{s!}\,(m+s)\ldots (m+1)\,\le\,
\frac{(m+s)^{s}}{s!}\,\le\,\frac{2^{s-1}(m^{s} +
s^{s})}{(s/e)^{s\mathstrut}}\,=\\
&
=\,\frac{1}{2}\Bigl(\Bigl(\frac{2em}{s}\Bigr)^{s}\,+\,(2e)^{s}\Bigr)\,\le\,(2em)^{s}\,\le\,\Bigl(\frac{2e\ln
x}{\ln 2}\Bigr)^{s}\,<\,(8\ln x)^{s}.
\end{align*}
Further,
\begin{align*}
& D_{2}(x)\,=\,\sum\limits_{k =
0}^{+\infty}\;\sum\limits_{2^{k}x<\delta\le
2^{k+1}x}\frac{1}{\delta}\;\le\;\sum\limits_{k =
0}^{+\infty}\frac{1}{2^{k}x}\sum\limits_{\delta \le
2^{k+1}x}1\,\le\, \sum\limits_{k =
0}^{+\infty}\frac{D_{1}(2^{k+1}x)}{2^{k}x}\,<\\
& <\,\frac{1}{x}\sum\limits_{k =
0}^{+\infty}\frac{(8\ln(2^{k+1}x))^{s}}{2^{k}}\,\ll\,\frac{1}{x}\sum\limits_{k
= 0}^{+\infty}\frac{(\ln
x)^{s}\,+\,(k+1)^{s}}{2^{k}}\,\ll_{s}\,\frac{(\ln x)^{s}}{x}.
\end{align*}
It remains to note that $s\le \tau(d)$.

\vspace{0.2cm}

\textbf{Lemma 12.} \emph{For any fixed} $a\ge 1$
\[
S_{a}(x)\,=\,C(\ln x)\beta(a)E_{a}(x)\,+\,\theta
R_{a}(x)\,+\,O(x\ln\ln x),
\]
\emph{where the constant} $C$ \emph{is defined in lemma 10},
\begin{align*}
& E_{a}(x)\,=\,\sum\limits_{n\le x}\frac{e_{a}(x)}{\tau(n)},\quad
e_{a}(n)\,=\,\frac{1}{\beta(a)}\sum\limits_{d|(a,n)}\beta\Bigl(\frac{an}{d^{2\mathstrut}}\Bigr),\\
& R_{a}(x)\,=\,\sum\limits_{d|a}\sum\limits_{\delta\le (\ln
x)^{3\mathstrut}}R_{a,d,\delta}(x),\\
& R_{a,d,\delta}(x)\,=\,\sum\limits_{\substack{q\le y/d \\ (q,
\frac{\Z a\delta}{\Z d}) =
1}}\frac{1}{\vf(q)}\sum\limits_{\chi\ne\chi_{0}}\Bigl|\sum\limits_{\substack{m\le
\frac{\Z x}{\Z d\delta\mathstrut} \\
(m,d)=1}}\frac{\chi(m)}{\tau(m)}\Bigr|,
\end{align*}
$y = \sqrt{x}(\ln x)^{-A}$, $A>0$ \emph{is arbitrary fixed number
and} $\delta$ \emph{runs through the sequence defined in lemma 11}.

\vspace{0.2cm}

\emph{Proof.} By definition of $\tau(n)$,
\begin{align*}
& S_{a}(x)\,=\,\sum\limits_{uv\le
x}\frac{1}{\tau(uv+a)}\,=\,\Bigl(\sum\limits_{u\le
\sqrt{x}}\;\;\sum\limits_{v\le\frac{\Z x}{\Z
u}}\,+\,\sum\limits_{\sqrt{x}<u\le x}\;\;\sum\limits_{v\le \frac{\Z
x}{\Z u}}\Bigr)\frac{1}{\tau(uv+a)}\,=\\
& =\,\Bigl(\sum\limits_{u\le\sqrt{x}}\sum\limits_{v\le \frac{\Z
x}{\Z u}}\,+\,\sum\limits_{\sqrt{x}<u\le
x}\;\;\sum\limits_{v\le \frac{\Z x}{\Z u}}\;\Bigr)\frac{1}{\tau(uv+a)}\,=\\
& =\,\Bigl(\sum\limits_{u\le\sqrt{x}}\sum\limits_{v\le \frac{\Z
x}{\Z u}}\,+\,\sum\limits_{v\le\sqrt{x}}\sum\limits_{u\le \frac{\Z
x}{\Z v}}\,-\,\sum\limits_{u\le\sqrt{x}}\sum\limits_{v\le
\sqrt{x}}\;\Bigr)\frac{1}{\tau(uv+a)}\,=\\
& =\,2\sum\limits_{u\le\sqrt{x}}\sum\limits_{v\le \frac{\Z x}{\Z
u}}\frac{1}{\tau(uv+a)}\,-\,\theta\,\tfrac{1}{2}\,x\,=\,2\sum\limits_{q\le\sqrt{x}}\sum\limits_{\substack{a<n\le
x+a \\ n\equiv a (\mmod
q)}}\frac{1}{\tau(n)}\,-\,\theta\,\tfrac{1}{2}\,x.
\end{align*}
Replacing the domain of $n$ in the inner sum to the interval $1\le
n\le x$, we obtain
\begin{align*}
&
S_{a}(x)\,=\,2\sum\limits_{q\le\sqrt{x}}\Bigl(\sum\limits_{\substack{1\le
n\le x \\ n\equiv a (\mmod q)}}\,-\,\sum\limits_{\substack{1\le n\le
a \\ n\equiv a (\mmod q)}}\,+\,\sum\limits_{\substack{x< n\le x+a \\
n\equiv a (\mmod
q)}}\Bigr)\frac{1}{\tau(n)}\,-\,\theta\,\tfrac{1}{2}\,x\,=\\
& =\,2\sum\limits_{q\le\sqrt{x}}\sum\limits_{\substack{1\le n\le x
\\ n\equiv a (\mmod
q)}}\frac{1}{\tau(n)}\,+\,2\theta_{1}\sum\limits_{q\le
\sqrt{x}}\Bigl(\frac{a}{q}+1\Bigr)\,-\,\theta\,\tfrac{1}{2}\,x\,=\\
& =\,2\sum\limits_{q\le\sqrt{x}}\sum\limits_{\substack{1\le n\le x
\\ n\equiv a (\mmod q)}}\frac{1}{\tau(n)}\,+\,\theta_{2}x.
\end{align*}

\vspace{0.2cm}

Suppose $A$ is an arbitrary positive number and let $y =
\sqrt{x}(\ln x)^{-A}$. Then the summands corresponding to the values
$y<q\le \sqrt{x}$ do not exceed
\begin{align*}
& 2\sum\limits_{y<q\le\sqrt{x}}\sum\limits_{\substack{n\le x \\
n\equiv a (\mmod
q)}}\frac{1}{2}\,\le\,\sum\limits_{y<q\le\sqrt{x}}\Bigl(\frac{x}{q}\,+\,1\Bigr)\,<\\
& <\,x\bigl(\ln\sqrt{x}\,-\,\ln
y\,+\,O(y^{-1})\bigr)\,+\,\sqrt{x}\,=\,A\ln\ln x\,+\,2\sqrt{x}(\ln
x)^{A}.
\end{align*}
Thus we get
\[
S_{a}(x)\,=\,V_{a}(x)\,+\,2\theta A\ln\ln x,\quad
V_{a}(x)\,=\,2\sum\limits_{q\le y}\sum\limits_{\substack{1\le n\le x
\\ n\equiv a (\mmod q)}}\frac{1}{\tau(n)}.
\]

Let us transform the last sum. Suppose $d = (q,a)$. Then $q =
dq_{1}$, $a = da_{1}$ where $(q_{1}, a_{1}) = 1$. If $n\equiv a
(\mmod q)$ then for some $k\ge 0$ we obtain
\[
n\,=\,a+kq\,=\,d(a_{1}\,+\,kq_{1})\,=\,dm,
\]
where $m\equiv a_{1} (\mmod q_{1})$ and $m \le x_{d} = x/d$. Thus
the inner sum takes the form
\[
\sum\limits_{\substack{m\le x_{d} \\ m\equiv a_{1} (\mmod
q_{1})}}\frac{1}{\tau(dm)}.
\]
Since all possible values of $d$ are among the divisors of $a$, it
follows that
\begin{align*}
& V_{a}(x)\,=\,2\sum\limits_{d|a}\sum\limits_{\substack{q\le y \\ (q,a) = d}}\;\sum\limits_{\substack{m\le x_{d} \\
m\equiv a_{1} (\mmod q_{1})}}\frac{1}{\tau(dm)}\,=\\
& =\,2\sum\limits_{d|a}\sum\limits_{\substack{q_{1}\le y_{d} \\
(q_{1},a_{1}) = 1}}\;\sum\limits_{\substack{n\le x_{d} \\
n\equiv a_{1} (\mmod q_{1})}}\frac{1}{\tau(dn)},
\end{align*}
where $y_{d} = y/d$.

Suppose $d$ is an arbitrary fixed divisor of $a$ and $\delta$ runs
through an increasing sequence that contains $1$ and positive
integers all whose prime divisors are among prime divisors of $d$.
(In particular, if $d = p$ is prime then $\delta$ take values $1, p,
p^{2}, p^{3}, \ldots$). Then for any integer $n$ there exists a
unique representation in the form $\delta m$ where $\delta$ belongs
to the above sequence and $(m,d) = (m,\delta) = 1$. Obviously for
such $n = \delta m$ we have $\tau(dn) = \tau(d\delta m) =
\tau(d\delta)\tau(m)$. Thus we obtain
\[
V_{a}(x)\,=\,2\sum\limits_{d|a}\sum\limits_{\delta}V_{a,d,\delta}(x)
\]
where
\[
V_{a,d,\delta}(x)\,=\,\sum\limits_{\substack{q_{1}\le y_{d} \\
(q_{1},a_{1}) = 1}}\;\;\sum\limits_{\substack{m\le \frac{\Z
x_{d}}{\Z \delta\mathstrut},\;(m,d)=1 \\ \delta m\equiv a_{1} (\mmod
q_{1})}}\frac{1}{\tau(d\delta m)}.
\]
Suppose that $(q_{1},a_{1})=1$ and the congruence $\delta m\equiv
a_{1} (\mmod q_{1})$ holds. If both the numbers $\delta$ and $q_{1}$
have the same divisor $\delta'>1$ then $\delta'$ divides $a_{1}$.
Therefore $(q_{1},a_{1})\ge \delta' >1$. This contradiction shows
that $(q_{1},\delta) = 1$. Thus the solutions of the above
congruence have the form $m\equiv a_{1}\delta^{*} (\mmod q_{1})$
This implies that
\[
V_{a,d,\delta}(x)\,=\,\sum\limits_{\substack{q_{1}\le y_{d} \\
(q_{1},a_{1}) = (q_{1},\delta) = 1}}\;\;\sum\limits_{\substack{m\le
\frac{\Z x_{d}}{\Z \delta\mathstrut},\,(m,d) = 1 \\ m\equiv
a_{1}\delta^{*} (\mmod q_{1})}}\frac{1}{\tau(d\delta m)}.
\]
Let $\delta_{0} = (\ln x)^{3}$. All the sums $V_{a,d,\delta}(x)$
corresponding to $\delta > \delta_{0}$ are estimated trivially:
\begin{align*}
& V_{a,d,\delta}(x)\,\le\, \sum\limits_{q_{1}\le
y_{d}}\sum\limits_{\substack{m\le \frac{\Z x_{d}}{\Z
\delta\mathstrut} \\ m\equiv a_{1}\delta^{*} (\mmod
q_{1})}}1\,\le\,\sum\limits_{q_{1}\le
y_{d}}\Bigl(\frac{x_{d}}{\delta
q_{1}}\,+\,1\Bigr)\,\ll\,\frac{x}{d\delta}\ln x\,+\,\frac{y}{d}
\end{align*}
where the constants in the symbols $\ll$ are absolute. By lemma 11,
the contribution of these terms to $V_{a}(x)$ do not exceed in order
\begin{align*}
&
\sum\limits_{d|a}\sum\limits_{\delta>\delta_{0}}\Bigl(\frac{x}{d\delta}\ln
x\,+\,\frac{y}{d}\Bigr)\,\ll\,x\ln
x\,\Bigl(\sum\limits_{d|a}\frac{1}{d}\Bigr)\sum\limits_{\delta>\delta_{0}}\frac{1}{\delta}\,+\,y\Bigl(\sum\limits_{d|a}\frac{1}{d}\Bigr)
\Bigl(\sum\limits_{\delta_{0}<\delta\le x}1\Bigr)\,\ll\\
& \ll\, x\ln x\,\frac{(\ln\delta_{0})^{\nu}}{\delta_{0}}\,+\,y(\ln
x)^{\nu}\,\ll\,\frac{x(\ln\ln x)^{\nu}}{(\ln
x)^{2}}\,\ll\,\frac{x}{\ln x},
\end{align*}
where $\nu = \tau(a)$ and the constants in the symbols $\ll$ depend
only on $a$. Thus,
\[
V_{a}(x)\,=\,2\sum\limits_{d|a}\sum\limits_{\delta\le\delta_{0}}V_{a,d,\delta}(x)\,+\,O\Bigl(\frac{x}{\ln
x}\Bigr).
\]
Next, for $1\le\delta\le \delta_{0}$ we have
\begin{align*}
& V_{a,d,\delta}(x)\,=\,\sum\limits_{\substack{q_{1}\le y_{d} \\
(q_{1}, a_{1}\delta) = 1}}\sum\limits_{\substack{m \le \frac{\Z
x_{d}}{\Z \delta\mathstrut} \\ (m,d) =
1}}\Bigl(\frac{1}{\vf(q_{1})}\sum\limits_{\chi \mmod
q_{1}}\overline{\chi}(a_{1}\delta^{*})\chi(m)\Bigr)\frac{1}{\tau(d\delta
m)}\,=\\
& =\,\sum\limits_{\substack{q_{1}\le y_{d} \\
(q_{1}, a_{1}\delta) = 1}}\frac{1}{\vf(q_{1})}\sum\limits_{\chi
\mmod
q_{1}}\chi(\delta)\overline{\chi}(a_{1})\sum\limits_{\substack{m\le
\frac{\Z x_{d}}{\Z \delta \mathstrut} \\ (m,d) =
1}}\frac{\chi(m)}{\tau(d\delta m)}\,=\\
& =\,\sum\limits_{\substack{q_{1}\le y_{d} \\
(q_{1}, a_{1}\delta) =
1}}\frac{1}{\vf(q_{1})}\sum\limits_{\substack{m\le \frac{\Z
x_{d}}{\Z \delta \mathstrut} \\ (m,d) = (m,q_{1}) =
1}}\frac{1}{\tau(d\delta m)}\;+\\
& +\,\frac{1}{\tau(d\delta)}\sum\limits_{\substack{q_{1}\le y_{d} \\
(q_{1}, a_{1}\delta) =
1}}\frac{1}{\vf(q_{1})}\sum\limits_{\substack{\chi \mmod q_{1}
\\ \chi\ne \chi_{0}}}\chi(\delta)\overline{\chi}(a_{1})\sum\limits_{\substack{m\le
\frac{\Z x_{d}}{\Z \delta \mathstrut} \\ (m,d) =
1}}\frac{\chi(m)}{\tau(m)}.
\end{align*}
Therefore,
\[
V_{a}(x)\,=\,W_{a}(x)\,+\,2\theta R_{a}(x)\,+\,O\Bigl(\frac{x}{\ln
x}\Bigr)
\]
where
\begin{align*}
& W_{a}(x)\,=\,2\sum\limits_{d|a}\sum\limits_{\delta\le\delta_{0}}\sum\limits_{\substack{q_{1}\le y_{d} \\
(q_{1}, a_{1}\delta) =
1}}\frac{1}{\vf(q_{1})}\sum\limits_{\substack{m\le \frac{\Z
x_{d}}{\Z \delta \mathstrut} \\ (m,d) = (m,q_{1}) =
1}}\frac{1}{\tau(d\delta m)},\\
&
R_{a}(x)\,=\,\sum\limits_{d|a}\sum\limits_{\delta\le\delta_{0}}R_{a,d,\delta}(x),\\
& R_{a,d,\delta}(x)\,=\,\sum\limits_{q\le
y_{d}}\frac{1}{\vf(q_{1})}\sum\limits_{\substack{\chi \mmod q
\\ \chi\ne \chi_{0}}}\Bigl|\sum\limits_{\substack{m\le \frac{\Z
x_{d}}{\Z \delta\mathstrut} \\ (m,d) =
1}}\frac{\chi(m)}{\tau(m)}\Bigr|.
\end{align*}
Let us replace the domain of $\delta$ in $W_{a}(x)$ by the segment
$1\le \delta\le x_{d}$. Additional summands do not exceed
\begin{align*}
& 2\sum\limits_{d|a}\sum\limits_{\delta_{0}<\delta\le
x_{d}}\sum\limits_{q_{1}\le
y_{d}}\frac{1}{\vf(q_{1})}\sum\limits_{m\le \frac{\Z x_{d}}{\Z
\delta\mathstrut}}\frac{1}{2}\,\ll\,\sum\limits_{d|a}\sum\limits_{\delta_{0}<\delta\le
x_{d}}\frac{x}{d\delta}\,\ln x\,\ll\\
& \ll\,x\ln
x\,\Bigl(\sum\limits_{d|a}\frac{1}{d}\Bigr)\sum\limits_{\delta>\delta_{0}}\frac{1}{\delta}\,\ll\,
x\ln x\,\frac{(\ln\delta_{0})^{\nu}}{\delta_{0}}\,\ll\,\frac{x}{\ln
x}.
\end{align*}
Thus,
\begin{align*}
& W_{a}(x)\,=\,2\sum\limits_{d|a}\sum\limits_{\delta}\sum\limits_{\substack{q_{1}\le y_{d} \\
(q_{1}, a_{1}\delta) =
1}}\frac{1}{\vf(q_{1})}\sum\limits_{\substack{m\le \frac{\Z
x_{d}}{\Z \delta \mathstrut} \\ (m,d) = (m,q_{1}) =
1}}\frac{1}{\tau(d\delta m)}\,+\,O\Bigl(\frac{x}{\ln x}\Bigr).
\end{align*}
Now our purpose is to transform double sum over $\delta$ and
$m\le\frac{\D x_{d}}{\D \delta\mathstrut}$ into single inner sum. We
put $n = \delta m$. Since for any $n$, $1\le n\le x_{d}$, there
exists a unique representation of the form $n = \delta m$ where
$(m,d)= 1$ and $\delta$ belongs to the above sequence, the condition
$(m,d) = 1$ in the inner sum in $W_{a}(x)$ may be omitted. Next, the
condition $(m,q_{1}) = 1$ should be replaced by the condition
$(q_{1},n) = 1$. Indeed, the equality $(q_{1},a_{1}\delta) = 1$ in
$W_{a}(x)$ implies that the sum over $\delta$ contains the terms
that obey the condition $(q_{1},\delta) = 1$. Then both the
conditions $(q_{1},m) = 1$ and $(q_{1},\delta) = 1$ are equivalent
to the single condition $(q_{1},n) = 1$. Thus we get
\begin{align*}
& W_{a}(x)\,=\,2\sum\limits_{d|a}\sum\limits_{\substack{q_{1}\le y_{d} \\
(q_{1}, a_{1}) = 1}}\frac{1}{\vf(q_{1})}\sum\limits_{\substack{n\le
x_{d} \\ (q_{1},n) = 1}}\frac{1}{\tau(dn)}\,+\,O\Bigl(\frac{x}{\ln
x}\Bigr)\,=\\
& =\,2\sum\limits_{d|a}\sum\limits_{n\le x_{d}}\frac{1}{\tau(dn)}\sum\limits_{\substack{q_{1}\le y_{d} \\
(q_{1}, na_{1}) = 1}}\frac{1}{\vf(q_{1})}\,+\,O\Bigl(\frac{x}{\ln
x}\Bigr).
\end{align*}
Applying the consequence of Lemma 9, we obtain
\[
W_{a}(x)\,=\,2\sum\limits_{d|a}\sum\limits_{n\le
x_{d}}\frac{1}{\tau(dn)}\bigl(C\beta(na_{1})\ln
y_{d}\,+\,O\bigl(\tau(na_{1})\bigr)\bigr)\,+\,O\Bigl(\frac{x}{\ln
x}\Bigr).
\]
Since $C\beta(na_{1})\ln d + O\bigl(\tau(na_{1})\bigr) \ll_{a}
\tau(n)$, it follows that
\begin{align*}
& W_{a}(x)\,=\\
& =\,2\sum\limits_{d|a}\sum\limits_{n\le
x_{d}}\frac{1}{\tau(dn)}\bigl(C\beta(na_{1})\bigl(\tfrac{1}{2}\ln
x\,-\,A\ln\ln
x\bigr)\,+\,O\bigl(\tau(n)\bigr)\bigr)\,+\,O\Bigl(\frac{x}{\ln
x}\Bigr)\,=\\
& =\,C(\ln x)\sum\limits_{d|a}\sum\limits_{n\le
x_{d}}\frac{\beta(na_{1})}{\tau(dn)}\,+\,r_{a}(x),
\end{align*}
where
\begin{align*}
& r_{a}(x)\,\ll\,(\ln\ln x)\sum\limits_{d|a}\sum\limits_{n\le
x_{d}}\frac{1}{\tau(dn)}\,+\,\sum\limits_{d|a}\sum\limits_{n\le
x_{d}}\frac{\tau(n)}{\tau(dn)}\,+\,\frac{x}{\ln x}\,\ll \\
& \ll\, (\ln\ln x)\Bigl(\sum\limits_{d|a}1\Bigr)\sum\limits_{m\le
x}\frac{1}{\tau(m)}\,+\,\sum\limits_{d|a}\frac{x}{d}\,+\,\frac{x}{\ln
x}\,\ll\\
& \ll\,\frac{x\ln\ln x}{\sqrt{\ln x}}\,+\,x\,+\,\frac{x}{\ln
x}\,\ll\,x.
\end{align*}
Thus,
\[
W_{a}(x)\,=\,C(\ln x)\sum\limits_{d|a}\sum\limits_{n\le
x/d}\frac{\beta\Bigl(\frac{\D na}{\D
d\mathstrut}\Bigr)}{\tau(nd)}\;+\;O(x).
\]
Changing the order of summation, we deduce that
\begin{align*}
& W_{a}(x)\,=\,C(\ln x)\sum\limits_{m\le
x}\frac{1}{\tau(m)}\sum\limits_{d|a,\;d|m}\beta\Bigl(\frac{\D na}{\D
d\mathstrut}\Bigr)\;+\;O(x)\,=\\
& =\,C(\ln x)\,\beta(a)\sum\limits_{m\le
x}\frac{e_{a}(m)}{\tau(m)}\;+\;O(x),
\end{align*}
where
\[
e_{a}(m)\,=\,\frac{1}{\beta(a)}\sum\limits_{d|(a,m)}\beta\Bigl(\frac{\D
na}{\D d^{2\mathstrut}}\Bigr)
\]
This completes the proof of the lemma.

\vspace{0.2cm}

The following lemma is the main assertion of the paper.

\vspace{0.2cm}

\textbf{Lemma 13.} \emph{Let} $d$ \emph{be a fixed integer}, $d\ge
1$. \emph{Then for any} $B > 0$ \emph{there exists} $A = A(B) > 0$
\emph{such that the estimation}
\[
R\,=\,\sum\limits_{q\le
Q}\frac{1}{\vf(q)}\;\sum\limits_{\substack{\chi \mmod q \\ \chi\ne
\chi_{0}}}\;\Bigl|\sum\limits_{\substack{n\le N \\ (n,d) =
1}}\frac{\chi(n)}{\tau(n)}\Bigr|\,\ll\, x(\ln x)^{-B}
\]
\emph{holds for any} $Q, N$ \emph{such that} $Q\le \sqrt{x}(\ln
x)^{-A}$, $N\le x$.

\vspace{0.2cm}

\textbf{Remark.} This assertion holds true for $A =
\tfrac{16}{3}\,B+4$.

\vspace{0.2cm}

\emph{Proof.} First we prove that if $T$ do not coincide with an
ordinate of a zero of $L(s,\chi)$ then the following inequality
holds:
\begin{align*}
& \Bigl|\sum\limits_{\substack{n\le N \\ (n,d) =
1}}\frac{\chi(n)}{\tau(n)}\Bigr|\,\ll\\
&\ll\,(\ln x)^{\alpha}\biggl(\frac{N}{T}\ln
x\,+\,\frac{N}{T}\int_{\sigma_{1}}^{\sigma_{2}}\bigl(|L(\sigma +
iT,\chi)|^{\frac{\Z 1}{\Z 2\mathstrut}}\,+\,|L(\sigma +
iT,\overline{\chi})|^{\frac{\Z 1}{\Z 2\mathstrut}}\bigr)d\sigma\;+\\
& +\,\sqrt{N}\int_{0}^{T}\bigl(|L(\sigma_{1} + it,\chi)|^{\frac{\Z
1}{\Z 2\mathstrut}}\,+\,|L(\sigma_{1} +
it,\overline{\chi})|^{\frac{\Z 1}{\Z
2\mathstrut}}\bigr)\frac{dt}{t+1}\;+\\
& +\,\sum\limits_{|\gamma|\le
T}\frac{1}{|\gamma|+1}\int_{\sigma_{1}}^{\beta}|L(\sigma +
i\gamma,\chi)|^{\frac{\Z 1}{\Z
2\mathstrut}}N^{\sigma}d\sigma\biggr).
\end{align*}
Here $\alpha = \tfrac{1}{24}$, $\sigma_{1} = \tfrac{1}{2}+(\ln
x)^{-1}$, $\sigma_{2} = \sigma_{1} + \tfrac{1}{2}$ and $\varrho =
\beta + i\gamma$ runs through all zeros of $L(s,\chi)$ in the
rectangle $\sigma_{1}<\beta\le 1$, $|\gamma|\le T$. Assume that
$\RRe s > 1$. Then the generating function $F(s;\chi,d) = F(s)$ of
the sequence $\chi(n)/\tau(n)$, $(n,d) = 1$, has the form
\[
F(s)\,=\,\sum\limits_{\substack{n = 1 \\ (n,d) =
1}}^{+\infty}\frac{\chi(n)}{\tau(n)}\,n^{-s}\,=\,\prod\limits_{p\nmid
d}F_{p}(s),\quad F_{p}(s)\,=\,1\,+\,\sum\limits_{k =
1}^{+\infty}\frac{\chi(p^{k})}{k+1}\,p^{-ks}.
\]
Taking for brevity $\kappa = \chi(p)p^{-s}$ we obtain for $\sigma
> 1$:
\[
F_{p}(s)\,=\,1\,+\,\tfrac{1}{2}\kappa\,+\,\tfrac{1}{3}\kappa^{2}\,+\,\tfrac{1}{4}\kappa^{3}\,+\ldots\,=\,(1-\kappa)^{-\frac{\Z
1}{\Z 2\mathstrut}}(1-\kappa^{2})^{\frac{\Z 1}{\Z
24\mathstrut}}(1-\kappa^{3})^{\frac{\Z 1}{\Z 24\mathstrut}}f_{p}(s)
\]
where the symbol $f_{p}(s)$ denotes a convergent series of the form
\[
1\,-\,\tfrac{49}{2880}\kappa^{4}\,-\,\tfrac{49}{1440}\kappa^{5}\,+\,\tfrac{1447}{362880}\kappa^{6}\,-\,\tfrac{3383}{120360}\kappa^{7}\,+\ldots\;.
\]
Therefore,
\begin{align*}
& F(s)\,=\,\prod\limits_{p\nmid
d}\Bigl(1\,-\,\frac{\chi(p)}{p^{s}}\Bigr)^{-\frac{\Z 1}{\Z
2\mathstrut}}\Bigl(1\,-\,\frac{\chi^{2}(p)}{p^{2s}}\Bigr)^{\frac{\Z
1}{\Z
24\mathstrut}}\Bigl(1\,-\,\frac{\chi^{3}(p)}{p^{3s}}\Bigr)^{\frac{\Z
1}{\Z 24\mathstrut}}f_{p}(s)\,=\\
& =\,\bigl(L(s,\chi)\bigr)^{\,\frac{\Z 1}{\Z
2\mathstrut}}\bigl(L(2s,\chi^{2})\bigr)^{-\,\frac{\Z 1}{\Z
24\mathstrut}}\bigl(L(3s,\chi^{3})\bigr)^{-\,\frac{\Z 1}{\Z
24\mathstrut}}\times\\
&
\times\,\prod\limits_{p|d}\Bigl(1\,-\,\frac{\chi(p)}{p^{s}}\Bigr)^{\frac{\Z
1}{\Z
2\mathstrut}}\Bigl(1\,-\,\frac{\chi^{2}(p)}{p^{2s}}\Bigr)^{-\frac{\Z
1}{\Z
24\mathstrut}}\Bigl(1\,-\,\frac{\chi^{3}(p)}{p^{3s}}\Bigr)^{-\frac{\Z
1}{\Z 24\mathstrut}}\prod\limits_{p\nmid d}f_{p}(s)\,=\\
& =\,\frac{\bigl(L(s,\chi)\bigr)^{\,\frac{\Z 1}{\Z
2\mathstrut}}}{\bigl(L(2s,\chi^{2})L(3s,\chi^{3})\bigr)^{\alpha}}\,\Phi(s),
\end{align*}
where $\alpha = \tfrac{1}{24}$ and the function $\Phi(s)$ is regular
in the half-plane $\RRe s > \tfrac{1}{4}$. Suppose $T\ge 2$ do not
coincide with an ordinate of a zero of $L(s,\chi)$. Taking $N_{1} =
[N]+\tfrac{1}{2}$, $\sigma_{2} = 1 + (\ln x)^{-1}$, by Perron's
formula we get
\[
\sum\limits_{\substack{n\le N \\ (n,d) =
1}}\frac{\chi(n)}{\tau(n)}\,=\,\frac{1}{2\pi
i}\int_{\sigma_{2}-iT}^{\sigma_{2}+iT}F(s)\frac{N_{1}^{s}}{s}\,ds\;+\;O\Bigl(\frac{N}{T}\ln
x\Bigr)\,+\,O(1).
\]
Suppose $\varrho = \beta + i\gamma$ runs through all the zeros of
$L(s,\chi)$ in the domain $\sigma_{2}<\beta\le 1$, $|\gamma|\le T$.
Let $\Gamma$ be a boundary of the rectangle with vertices
$\sigma_{1}\pm iT$, $\sigma_{2}\pm iT$ and with horizontal cuts
going from the left side of the rectangle to each zero $\varrho$.
Applying Cauchy's theorem we obtain
\[
\frac{1}{2\pi
i}\int_{\sigma_{2}-iT}^{\sigma_{2}+iT}F(s)\frac{N_{1}^{s}}{s}\,ds\;=\;-\Bigl(I_{1}\,+\,I_{2}\,+\,I_{3}\,+\,\sum\limits_{|\gamma|\le
T}I(\varrho)\Bigr)
\]
where the symbols $I_{1}, I_{2}, I_{3}$ denote the integrals
\[
\frac{1}{2\pi
i}\int_{\sigma_{2}+iT}^{\sigma_{1}+iT}F(s)\frac{N_{1}^{s}}{s}\,ds,\;\;
\frac{1}{2\pi
i}\int_{\sigma_{1}-iT}^{\sigma_{2}-iT}F(s)\frac{N_{1}^{s}}{s}\,ds,
\;\; \frac{1}{2\pi
i}\,\text{v.p.}\,\int_{\sigma_{1}+iT}^{\sigma_{1}-iT}F(s)\frac{N_{1}^{s}}{s}\,ds
\]
respectively, and $I(\varrho)$ means the sum of the integrals over
the upper and lower edges of the cut:
\[
I(\varrho)\,=\,\frac{1}{2\pi i}\biggl(\int_{\sigma_{1}+i(\gamma +
0)}^{\beta + i(\gamma + 0)}\;+\;\int_{\beta + i(\gamma -
0)}^{\sigma_{1}+i(\gamma - 0)}\biggr)F(s)\frac{N_{1}^{s}}{s}\,ds.
\]
Since the infinite product for $\Phi(s)$ converges absolutely in the
half-plane $\RRe s > \tfrac{1}{4}$, then $|\Phi(s)| = O(1)$ along
the contour $\Gamma$. Moreover, for $\sigma\ge \sigma_{1}$ we obtain
\begin{align*}
& |L(3s,\chi^{3})|^{-1}\,=\,\Bigl|\sum\limits_{n =
1}^{+\infty}\frac{\mu(n)\chi^{3}(n)}{n^{3s}}\Bigr|\,\le\,\sum\limits_{n
= 1}^{+\infty}n^{-3/2}\,=\,O(1),\\
& |L(2s,\chi^{2})|^{-1}\,=\,\Bigl|\sum\limits_{n =
1}^{+\infty}\frac{\mu(n)\chi^{2}(n)}{n^{2s}}\Bigr|\,\le\,\sum\limits_{n
= 1}^{+\infty}n^{-2\sigma}\,\le\,\sum\limits_{n =
1}^{+\infty}n^{2\sigma_{1}}\,\le\\
& \le\,
1\,+\,\int_{1}^{+\infty}u^{-2\sigma_{1}}du\,=\,\frac{2\sigma_{1}}{2\sigma_{1}-1}\,=\,\tfrac{1}{2}\,\ln
x + 1\,<\,\ln x.
\end{align*}
Thus,
\[
|F(s)|\,\ll\,(\ln x)^{\alpha}|L(s,\chi)|^{\frac{\Z 1}{\Z
2\mathstrut}}
\]
along the contour $\Gamma$. Therefore,
\begin{align*}
& I_{1}\ll\,\int_{\sigma_{1}}^{\sigma_{2}}|F(\sigma +
iT)|\frac{N_{1}^{\sigma}d\sigma}{\sqrt{\sigma^{2}+T^{2}}}\,\ll\,\frac{N^{\sigma_{2}}}{T}(\ln
x)^{\alpha}\int_{\sigma_{1}}^{\sigma_{2}}|L(\sigma +
iT,\chi)|^{\frac{\Z 1}{\Z 2\mathstrut}}d\sigma\,\ll\\
& \ll\,\frac{N}{T}(\ln
x)^{\alpha}\int_{\sigma_{1}}^{\sigma_{2}}|L(\sigma +
iT,\chi)|^{\frac{\Z 1}{\Z 2\mathstrut}}d\sigma,\\
& I_{2}\,\ll\,\frac{N}{T}(\ln
x)^{\alpha}\int_{\sigma_{1}}^{\sigma_{2}}|L(\sigma +
iT,\overline{\chi})|^{\frac{\Z 1}{\Z 2\mathstrut}}d\sigma.
\end{align*}
Next,
\begin{align*}
&
|I_{3}|\,=\,\frac{1}{2\pi}\Bigl|\text{v.p.}\int_{-T}^{T}F(\sigma_{1}+it
)\,\frac{N_{1}^{\sigma_{1}+it}}{\sigma_{1}+it}\,dt\Bigr|\,\ll\,N^{\sigma_{1}}\int_{-T}^{T}\frac{|F(\sigma_{1}+it)|dt}{\sqrt{\sigma_{1}^{2}+t^{2}}}\,\ll\\
& \ll\,\sqrt{N}(\ln
x)^{\alpha}\int_{-T}^{T}\frac{|L(\sigma_{1}+it,\chi)|^{\frac{\Z
1}{\Z 2\mathstrut}}dt}{\sqrt{\sigma_{1}^{2}+t^{2}}}\,\ll \\
& \ll\,\sqrt{N}(\ln
x)^{\alpha}\int_{0}^{T}\bigl(L(\sigma_{1}+it,\chi)|^{\frac{\Z 1}{\Z
2\mathstrut}}\,+\,L(\sigma_{1}+it,\overline{\chi})|^{\frac{\Z 1}{\Z
2\mathstrut}}\bigr)\frac{dt}{t+1}.
\end{align*}
Finally, each of the integrals $I(\varrho)$ obeys the inequality
\[
|I(\varrho)|\,\ll\,\frac{(\ln x)^{\alpha +
1}}{|\gamma|+1}\int_{\sigma_{1}}^{\beta}|L(\sigma_{1}+i\gamma,\chi)|^{\frac{\Z
1}{\Z 2\mathstrut}}N^{\sigma}d\sigma.
\]
Summing the above estimates we obtain the required inequality for
the sum of $\chi(n)/\tau(n)$, $n\le N$, $(n,d) = 1$. Summing these
inequalities over all non-principle characters $\chi \mmod q$ and
over $q\le Q$, we get
\[
R\,\ll\,(\ln x)^{\alpha}\sum\limits_{j = 1}^{4}R_{j},
\]
where
\begin{align*}
& R_{1}\,=\,\sum\limits_{q\le
Q}\frac{1}{\vf(q)}\sum\limits_{\chi\ne\chi_{0}}\frac{N}{T}\ln x,\\
& R_{2}\,=\,\sum\limits_{q\le
Q}\frac{1}{\vf(q)}\sum\limits_{\chi\ne\chi_{0}}\frac{N}{T}\int_{\sigma_{1}}^{\sigma_{2}}\bigl(|L(\sigma+iT,\chi)|^{\frac{\Z
1}{\Z 2\mathstrut}}\,+\,|L(\sigma+iT,\overline{\chi})|^{\frac{\Z
1}{\Z 2\mathstrut}}\bigr)d\sigma,\\
& R_{3}\,=\,\sum\limits_{q\le
Q}\frac{1}{\vf(q)}\sum\limits_{\chi\ne\chi_{0}}\sqrt{N}\int_{0}^{T}\bigl(|L(\sigma_{1}+it,\chi)|^{\frac{\Z
1}{\Z 2\mathstrut}}\,+\,|L(\sigma_{1}+it,\overline{\chi})|^{\frac{\Z
1}{\Z 2\mathstrut}}\bigr)\frac{dt}{t+1},\\
& R_{4}\,=\,\sum\limits_{q\le
Q}\frac{1}{\vf(q)}\sum\limits_{\chi\ne\chi_{0}}\sum\limits_{|\gamma|\le
T}\frac{1}{|\gamma|+1}\int_{\sigma_{1}}^{\beta}N^{\sigma}|L(\sigma+i\gamma,\overline{\chi})|^{\frac{\Z
1}{\Z 2\mathstrut}}d\sigma.
\end{align*}
The values of $T$ depends on $q$ and $\chi$ only: $T = T(q,\chi)$.
We choose $T(q,\chi)$ in a following way. The points $Q\cdot
2^{-k}$, $k = 0,1,2,\ldots$ split the segment $(1,Q]$ into the
intervals of the type $(M, M_{1}]$, where $M_{1}\le 2M$, $M_{1}\le
Q$. Suppose the character $\chi$ modulo $q$ is induced by a
primitive character $\chi_{1}$ modulo $q_{1}$, $q = q_{1}r$, where
$M<q_{1}\le M_{1}$. Then we put $T = M(\ln x)^{C}$, where the
constant $C$, $0<C<A$, will be chosen later. Therefore, it follows
that $T(q,\chi) = T(q_{1},\chi_{1})$; in particular, $T(q,\chi)$ do
not depends on $r$. Replacing (if it is necessary) the value
$T(q,\chi)$ by the value $T(q,\chi)+h$ for some $h$, $0<h\le c(\ln
x)^{-1}$, we may assume that $T(q,\chi)$ does not coincide with an
ordinate of a zero of $L(s,\chi)$. In each case, we obviously have
\[
M(\ln x)^{C} \le T(q,\chi)\le  M(\ln x)^{C} + c(\ln x)^{-1}\,\le\,
2M(\ln x)^{C}.
\]
In the following, the sums over $\chi \mmod q$, $\chi\ne\chi_{0}$
are replaced by the sums over primitive characters $\chi_{1} \mmod
q_{1}$.

\vspace{0.2cm}

$1^{\circ}$. Estimation of $R_{1}$. Obviously we have
\begin{align*}
R_{1}\,&\ll\,N(\ln x)\sum\limits_{q\le
Q}\frac{1}{\vf(q)}\sum\limits_{\chi\ne\chi_{0}}\frac{1}{T(q,\chi)}\,\ll\\
& \ll\,N(\ln x)\sum\limits_{q_{1}\le Q}\sum\limits_{r\le \frac{\Z
Q}{\Z q_{1}}}\frac{1}{\vf(q_{1}r)}\stsum\limits_{\chi_{1} \mmod
q_{1}}\frac{1}{T(q_{1},\chi_{1})}.
\end{align*}
Using the symbol $\prsum_{M\le Q}$ for the summation over the points
$M = Q\cdot 2^{-k}$, by the inequality $\vf(q_{1}r)\ge
\vf(q_{1})\vf(r)$ we obtain
\begin{align*}
& R_{1}\,\ll\,N(\ln x)\prsum\limits_{M\le
Q}\;\sum\limits_{M<q_{1}\le
M_{1}}\frac{1}{\vf(q_{1})}\sum\limits_{r\le \frac{\Z Q}{\Z
q_{1}}}\frac{1}{\vf(r)}\stsum\limits_{\chi_{1} \mmod
q_{1}}\frac{1}{M(\ln x)^{C}}\,\ll\\
& \ll\,\frac{N\ln x}{(\ln x)^{C}}\prsum\limits_{M\le
Q}\frac{1}{M}\sum\limits_{M<q_{1}\le
M_{1}}\frac{1}{\vf(q_{1})}\Bigl(\stsum\limits_{\chi_{1} \mmod
q_{1}}1\Bigr)\sum\limits_{r\le Q}\frac{1}{\vf(r)}\,\ll \\
& \ll\,\frac{N(\ln x)^{2}}{(\ln x)^{C}}\prsum\limits_{M\le
Q}\frac{1}{M}\Bigl(\sum\limits_{M<q_{1}\le
M_{1}}1\Bigr)\,\ll\,\frac{N(\ln x)^{2}}{(\ln
x)^{C}}\prsum\limits_{M\le Q}1\,\ll\,\frac{N(\ln x)^{3}}{(\ln
x)^{C}}.
\end{align*}

\vspace{0.2cm}

$2^{\circ}$. Estimation of $R_{2}$. Since $\overline{\chi}$ and
$\chi$ run through the same set of characters, we get
\[
R_{2}\,\ll\,N\int_{\sigma_{1}}^{\sigma_{2}}r_{2}(\sigma)d\sigma,
\]
where
\[
r_{2}(\sigma)\,=\,\sum\limits_{q\le
Q}\frac{1}{\vf(q)}\sum\limits_{\substack{\chi \mmod q \\ \chi\ne
\chi_{0}}}\frac{|L(\sigma+iT,\chi)|^{\frac{\Z 1}{\Z
2\mathstrut}}}{T(q,\chi)}.
\]
Using the same arguments as above, by lemma 7 we obtain ($T =
T(q_{1},\chi_{1})$):
\begin{align*}
& r_{2}(\sigma)\,\le\,\sum\limits_{q_{1}\le
Q}\frac{1}{\vf(q_{1})}\sum\limits_{r\le \frac{\Z Q}{\Z
q_{1}}}\frac{1}{\vf(r)}\stsum\limits_{\chi_{1} \mmod
q_{1}}\frac{\sqrt{\tau(r)}}{T(q_{1},\chi_{1})}|L(\sigma +
iT,\chi_{1})|^{\frac{\Z 1}{\Z 2\mathstrut}}\,\le\\
& \le\, \sum\limits_{q_{1}\le
Q}\frac{1}{\vf(q_{1})}\sum\limits_{r\le \frac{\Z Q}{\Z
q_{1}}}\frac{\sqrt{\tau(r)}}{\vf(r)}\stsum\limits_{\chi_{1} \mmod
q_{1}}\frac{|L(\sigma + iT,\chi_{1})|^{\frac{\Z
1}{\Z 2\mathstrut}}}{T(q_{1},\chi_{1})}\,\ll\\
& \ll\, \prsum\limits_{M\le Q}\;\sum\limits_{M<q_{1}\le
M_{1}}\frac{1}{\vf(q_{1})}\sum\limits_{r\le \frac{\Z Q}{\Z
q_{1}}}\frac{\sqrt{\tau(r)}}{\vf(r)}\stsum\limits_{\chi_{1} \mmod
q_{1}}\frac{|L(\sigma + iT,\chi_{1})|^{\frac{\Z 1}{\Z
2\mathstrut}}}{M(\ln x)^{C}}\,\ll\\
& \ll\,(\ln x)^{-C}\prsum\limits_{M\le
Q}\frac{1}{M}\sum\limits_{M<q_{1}\le
M_{1}}\frac{1}{\vf(q_{1})}\sum\limits_{r\le \frac{\Z Q}{\Z
q_{1}}}\frac{\sqrt{\tau(r)}}{\vf(r)}\stsum\limits_{\chi_{1} \mmod
q_{1}}|L(\sigma + iT,\chi_{1})|^{\frac{\Z 1}{\Z 2\mathstrut}}.
\end{align*}
By  the inequalities
\[
\sum\limits_{r\le X}\sqrt{\tau(r)}\le \Bigl(X\sum\limits_{r\le
X}\tau(r)\Bigr)^{\frac{\Z 1}{\Z 2\mathstrut}}\,\ll\,X\sqrt{\ln
X},\quad \frac{1}{\vf(r)}\,\ll\,\frac{\ln\ln(r+3)}{r},
\]
we deduce that
\begin{align*}
& \sum\limits_{r\le \frac{\Z Q}{\Z
q_{1}}}\frac{\sqrt{\tau(r)}}{\vf(r)}\,\ll\,(\ln\ln
Q)\sum\limits_{r\le Q}\frac{\sqrt{\tau(r)}}{r}\,\ll\\
& \ll\,(\ln\ln x)\Bigl(\frac{1}{Q}\sum\limits_{r\le
Q}\sqrt{\tau(r)}\,+\,\int_{2}^{Q}\Bigl(\sum\limits_{r\le
u}\sqrt{\tau(r)}\Bigr)\frac{du}{u^{2}}\,\Bigr)\,\ll\,(\ln\ln
x)\int_{2}^{Q}\!\sqrt{\ln u}\;\frac{du}{u}\\
& \ll\, (\ln x)^{\frac{\Z 3}{\Z 2\mathstrut}}\ln\ln x
\end{align*}
and therefore
\begin{align*}
& r_{2}(\sigma)\,\ll\,\frac{(\ln x)^{\frac{\Z 3}{\Z
2\mathstrut}}\ln\ln x}{(\ln x)^{C}}\prsum\limits_{M\le
Q}\frac{1}{M}\sum\limits_{M<q_{1}\le
M_{1}}\frac{1}{\vf(q_{1})}\stsum\limits_{\chi_{1} \mmod
q_{1}}|L(\sigma + iT,\chi_{1})|^{\frac{\Z 1}{\Z 2\mathstrut}}\,\ll
\\
& \ll\, \frac{(\ln x)^{\frac{\Z 3}{\Z 2\mathstrut}}(\ln\ln
x)^{2}}{(\ln x)^{C}}\prsum\limits_{M\le
Q}\frac{1}{M^{2}}\,r_{2}(\sigma, M),
\end{align*}
where
\[
r_{2}(\sigma, M)\,=\,\sum\limits_{M<q_{1}\le
M_{1}}\stsum\limits_{\chi_{1} \mmod q_{1}}|L(\sigma +
iT,\chi_{1})|^{\frac{\Z 1}{\Z 2\mathstrut}}.
\]
Applying H\"{o}lder's inequality, we obtain
\begin{align*}
& r_{2}^{4}(\sigma, M)\,\ll\, \Bigl(\sum\limits_{M<q_{1}\le
M_{1}}\stsum\limits_{\chi_{1} \mmod
q_{1}}1\Bigr)^{3}\sum\limits_{M<q_{1}\le
M_{1}}\stsum\limits_{\chi_{1} \mmod q_{1}}|L(\sigma +
iT,\chi_{1})|^{2}\,\ll\\
& \ll\,M^{6}\sum\limits_{M<q_{1}\le M_{1}}\stsum\limits_{\chi_{1}
\mmod q_{1}}|L(\sigma + iT,\chi_{1})|^{2}.
\end{align*}
Let $Y = 2M_{1}T(q_{1},\chi_{1})$; then for any $q_{1}$ such that
$M<q_{1}\le M_{1}$, the following inequality holds: $Y\ge
q_{1}(T(q_{1},\chi_{1})+1)/\pi$. By lemma 3, we get
\[
|L(\sigma + iT,\chi_{1})|^{2}\,\ll\,\Bigl|\sum\limits_{n\le
Y}\frac{\chi(n)}{n^{\sigma +
iT}}\Bigr|^{2}\,+\,(q_{1}Y^{-\sigma})^{2}.
\]
Since
\begin{align*}
&
q_{1}Y^{-\sigma}\,\le\,\frac{M_{1}}{\sqrt{Y}}\,=\,\frac{M_{1}}{\sqrt{2M_{1}T(q_{1},\chi_{1})}}
\,<\,\sqrt{\frac{M_{1}}{T(q_{1},\chi_{1})}}\,<\\
& <\,\sqrt{\frac{2M}{M(\ln x)^{C}}}\,=\,\sqrt{2}(\ln
x)^{-C/2}\,<\,1,
\end{align*}
for $\sigma \ge \tfrac{1}{2}$, we get
\begin{align*}
& r_{2}^{4}(\sigma, M)\,\ll\, M^{8}\,+\,M^{6}\sum\limits_{q_{1}\le
M_{1}}\stsum\limits_{\chi_{1} \mmod q_{1}}\Bigl|\sum\limits_{n\le
Y}\frac{\chi(n)}{n^{\sigma + iT}}\Bigr|^{2}.
\end{align*}
Setting $Q = M_{1}$, $M = 0$, $N = Y$ in lemma 2, we obtain
\[
r_{2}^{4}(\sigma,
M)\,\ll\,M^{8}\,+\,M^{6}(M^{2}\,+\,Y)\sum\limits_{n\le
Y}\frac{1}{n^{2\sigma}}\,\ll\,M^{8}(\ln x)^{C+1}
\]
and therefore
\begin{align*}
& r_{2}(\sigma)\,\ll\, \frac{(\ln x)^{\frac{\Z 3}{\Z
2\mathstrut}}(\ln\ln x)^{2}}{(\ln x)^{C}}\prsum\limits_{M\le
Q}\frac{1}{M^{2\mathstrut}}\,M^{2}(\ln x)^{\frac{\Z 1}{\Z
4\mathstrut}(C+1)}\,\ll\\
& \ll\, \frac{(\ln x)^{\frac{\Z 7}{\Z 4\mathstrut}}(\ln\ln
x)^{2}}{(\ln x)^{\frac{\Z 3}{\Z 4\mathstrut}C}}\prsum\limits_{M\le
Q}1\,\ll \frac{(\ln x)^{\frac{\Z 11}{\Z 4\mathstrut}}(\ln\ln
x)^{2}}{(\ln x)^{\frac{\Z 3}{\Z 4\mathstrut}C}}.
\end{align*}
Thus,
\[
R_{2}\,\ll\,N\,\frac{(\ln x)^{\frac{\Z 11}{\Z 4\mathstrut}}}{(\ln
x)^{\frac{\Z 3}{\Z 4\mathstrut}C}}\,(\ln\ln x)^{2}.
\]

\vspace{0.2cm}

$3^{\circ}$. Estimation of $R_{3}$. By the same arguments we obtain
\begin{align*}
& R_{3}\,\ll\,\sqrt{N}\sum\limits_{q\le
Q}\frac{1}{\vf(q)}\sum\limits_{\substack{\chi \mmod q \\ \chi\ne
\chi_{0}}}\int_{0}^{T(q,\chi)}|L(\sigma_{1}+it,\chi)|^{\frac{\Z
1}{\Z 2\mathstrut}}\,\frac{dt}{t+1}\,\ll \\
& \ll\,\sqrt{N}\sum\limits_{q_{1}\le
Q}\frac{1}{\vf(q_{1})}\sum\limits_{r\le\frac{\Z Q}{\Z
q_{1}}}\frac{\sqrt{\tau(r)}}{\vf(r)}\stsum\limits_{\chi_{1} \mmod
q_{1}}\int_{0}^{T(q_{1},\chi_{1})}|L(\sigma_{1}+it,\chi)|^{\frac{\Z
1}{\Z 2\mathstrut}}\,\frac{dt}{t+1}.
\end{align*}
Since the terms in the inner sum over $\chi_{1}$ do not depend on
$r$, we easily get
\begin{align*}
& R_{3}\,\ll\,\sqrt{N}(\ln x)^{\frac{\Z 3}{\Z 2\mathstrut}}(\ln\ln
x)\sum\limits_{q_{1}\le
Q}\frac{1}{\vf(q_{1})}\stsum\limits_{\chi_{1} \mmod
q_{1}}\int_{0}^{T(q_{1},\chi_{1})}|L(\sigma_{1}+it,\chi)|^{\frac{\Z
1}{\Z 2\mathstrut}}\,\frac{dt}{t+1}\,\ll\\
& \ll\,\sqrt{N}(\ln x)^{\frac{\Z 3}{\Z 2\mathstrut}}(\ln\ln
x)\prsum\limits_{M\le Q}\sum\limits_{M<q_{1}\le
M_{1}}\frac{1}{\vf(q_{1})}\stsum\limits_{\chi_{1} \mmod
q_{1}}\int_{0}^{T(q_{1},\chi_{1})}\frac{|L(\sigma_{1}+it,\chi)|^{\frac{\Z
1}{\Z 2\mathstrut}}dt}{t+1}\\
& \ll\,\sqrt{N}(\ln x)^{\frac{\Z 3}{\Z 2\mathstrut}}(\ln\ln
x)^{2}\prsum\limits_{M\le
Q}\frac{1}{M}\,\int_{0}^{T_{1}}\frac{r_{3}(t,M)}{t+1}\,dt,
\end{align*}
where $T_{1} = 2M(\ln x)^{C}$,
\[
r_{3}(t,M)\,=\,\sum\limits_{M<q_{1}\le M_{1}}\stsum\limits_{\chi_{1}
\mmod q_{1}}|L(\sigma_{1}+it,\chi_{1})|^{\frac{\Z 1}{\Z
2\mathstrut}}.
\]
From H\"{o}lder's inequality and lemma 2 it follows that
\[
r_{3}(t,M)\,\ll\,M^{2}(\ln x)^{\frac{\Z 1}{\Z 4\mathstrut}(C+1)}.
\]
Thus we find
\begin{align*}
& R_{3}\,\ll\,\sqrt{N}(\ln x)^{\frac{\Z 3}{\Z 2\mathstrut}}(\ln\ln
x)^{2}\prsum\limits_{M\le Q}\frac{1}{M}\, M^{2}(\ln x)^{\frac{\Z
1}{\Z 4\mathstrut}(C+1)}\,\ll\\
& \ll\,Q\sqrt{N}(\ln x)^{\frac{\Z 1}{\Z 4}(C+7)}(\ln\ln x)^{2}.
\end{align*}

\vspace{0.2cm}

$4^{\circ}$. Estimation of $R_{4}$. For given character $\chi \mmod
q$ and a zero $\varrho = \beta + i\gamma$ of $L(s,\chi)$ we define
the function $g_{\chi}(\varrho,\sigma)$ as follows:
\begin{equation*}
g_{\chi}(\varrho,\sigma)\,=\,
\begin{cases}
1,\quad & \text{if}\; \beta > \sigma,\\
0,\quad & \text{otherwise}.
\end{cases}
\end{equation*}
It's obvious that if $\chi_{1} \mmod q_{1}$ induces $\chi\mmod q$
then $g_{\chi_{1}}(\varrho,\sigma)$ coincides with
$g_{\chi}(\varrho,\sigma)$ for $\sigma\ge \tfrac{1}{2}$. Thus we
have
\begin{align*}
& R_{4}\,\le\,\sum\limits_{q\le
Q}\frac{1}{\vf(q)}\sum\limits_{\substack{\chi\mmod q \\
\chi\ne\chi_{0}}}\;\;\sum\limits_{\substack{|\gamma|\le T(q,\chi) \\
\beta >
0.5}}\frac{1}{|\gamma|+1}\int_{0.5}^{1}g_{\chi}(\varrho,\sigma)N^{\sigma}|L(\sigma
+ i\gamma,\chi)|^{\frac{\Z 1}{\Z 2\mathstrut}}\,d\sigma\,\le\\
& \le \int_{0.5}^{1}N^{\sigma}\sum\limits_{q_{1}\le
Q}\frac{1}{\vf(q_{1})}\sum\limits_{r\le\frac{\Z Q}{\Z
q_{1}}}\frac{\sqrt{\tau(r)}}{\vf(r)}\!\stsum\limits_{\chi_{1} \mmod
q_{1}}\;\sum\limits_{\substack{|\gamma|\le T(q,\chi) \\
\beta >
0.5}}\frac{g_{\chi_{1}}(\varrho,\sigma)}{|\gamma|+1}|L(\sigma +
i\gamma,\chi_{1})|^{\frac{\Z 1}{\Z 2\mathstrut}}\,d\sigma
\end{align*}
\begin{align*}
& \ll\,(\ln x)^{\frac{\Z 3}{\Z 2\mathstrut}}(\ln\ln
x)^{2}\!\int_{0.5}^{1}\!N^{\sigma}\!\sum\limits_{q_{1}\le
Q}\frac{1}{q_{1}}\!\stsum\limits_{\chi_{1} \mmod q_{1}}\;\sum\limits_{\substack{|\gamma|\le T(q,\chi) \\
\beta >
0.5}}\!\frac{g_{\chi_{1}}(\varrho,\sigma)}{|\gamma|+1}|L(\sigma +
i\gamma,\chi_{1})|^{\frac{\Z 1}{\Z 2\mathstrut}}\,d\sigma\\
& \ll\, (\ln x)^{\frac{\Z 3}{\Z 2\mathstrut}}(\ln\ln
x)^{2}\prsum\limits_{M\le Q}j(M),
\end{align*}
where
\begin{align*}
&
j(M)\,=\,\frac{1}{M}\int_{0.5}^{1}N^{\sigma}r_{4}(\sigma,M)d\sigma,\\
& r_{4}(\sigma, M)\,=\,\sum\limits_{M<q_{1}\le
M_{1}}\stsum\limits_{\chi_{1} \mmod q_{1}}\;\sum\limits_{\substack{|\gamma|\le T(q,\chi) \\
\beta >
0.5}}\!\frac{g_{\chi_{1}}(\varrho,\sigma)}{|\gamma|+1}|L(\sigma +
i\gamma,\chi_{1})|^{\frac{\Z 1}{\Z 2\mathstrut}}.
\end{align*}

In the below, we consider the cases of <<small>> and <<large>> $M$
separately.

\vspace{0.2cm}

$4\text{a}^{\circ}$. The case of <<small>> $M$: $1\le M\le (\ln
x)^{A}$. Obviously we have
\[
T(q_{1},\chi_{1}) \le 2M(\ln x)^{C}\le 2(\ln x)^{A+C}.
\]
If
$\varrho = \beta + i\gamma$ is a zero of $L(s,\chi_{1})$ and
$0<|\gamma|\le T(q_{1},\chi_{1})$ then lemma 4 implies that
\[
\beta\,\le\,1\,-\,\frac{c_{1}}{\ln{q_{1}(|\gamma|+2)}}\,\le\,1\,-\,\frac{c}{\ln{\bigl(2(\ln
c)^{2A+C}\bigr)}}\,<\,1\,-\,\frac{1}{\sqrt{\ln x}}.
\]
Further, if there exists a real zero $\beta$ of $L(s,\chi_{1})$ then
lemma 5 implies (with $\vep = (2A)^{-1}$) that
\[
\beta\,\le\,1\,-\,\frac{c}{q_{1}^{\vep}}\,\le\,1\,-\,\frac{c}{(\ln
x)^{A\vep}}\,=\,1\,-\,\frac{c}{\sqrt{\ln x}}
\]
for some $c > 0$. Without loss of generality, we may assume that
$0<c<1$. Then it follows that for any zero $\varrho$ of the function
$L(s,\chi_{1})$ under considering the inequality $|\gamma|\le
T(q_{1},\chi_{1})$ implies that $g_{\chi_{1}}(\varrho,\sigma) = 0$
for any $\sigma \ge 1 - \frac{\D 1}{\D \sqrt{\ln x}}$.

Now, if $\tfrac{1}{2}\le\sigma\le 1$ then from lemma 3 it follows
that
\begin{align*}
& |L(\sigma+it,\chi_{1})|\,\ll\,\Bigl|\sum\limits_{n\le
q_{1}(|t|+1)}\chi(n)n^{-\sigma-it}\Bigr|\,+\,q_{1}\bigl(q_{1}(|t|+1)\bigr)^{-\sigma}\,\\
& \ll\,\bigl(q_{1}(|t|+1)\bigr)^{1-\sigma}\ln{q_{1}(|t|+1)}.
\end{align*}
By lemma 7, we easily get
\begin{align*}
& |L(\sigma+i\gamma,\chi_{1})|^{\frac{\Z 1}{\Z
2\mathstrut}}\,\ll\,\bigl(q_{1}T(q_{1},\chi_{1})\bigr)^{\frac{\Z
1-\sigma}{\Z 2\mathstrut}}\ln\ln x\,\ll\,\bigl(M^{2}(\ln
x)^{C}\bigr)^{\frac{\Z 1-\sigma}{\Z 2\mathstrut}}\ln\ln x\,<\\
&<\, (\ln x)^{\frac{\Z 3}{\Z 4\mathstrut}A}\ln\ln x,\\
& r_{4}(\sigma,M)\,\ll\,(\ln x)^{\frac{\Z 3}{\Z 4\mathstrut}A}\ln\ln
x\,\sum\limits_{M<q\le M_{1}}\stsum\limits_{\chi_{1} \mmod
q_{1}}\;\sum\limits_{|\gamma|\le
T(q_{1},\chi_{1})}\frac{1}{|\gamma|+1}\,\ll\\
& <\,(\ln x)^{\frac{\Z 3}{\Z 4\mathstrut}A}\ln\ln
x\,\sum\limits_{M<q\le M_{1}}\stsum\limits_{\chi_{1} \mmod
q_{1}}\ln^{2}{q_{1}T(q_{1},\chi_{1})}\,\ll\,(\ln x)^{\frac{\Z 3}{\Z
4\mathstrut}A}(\ln\ln x)^{3}M^{2}\,\ll\\
& \ll\, M(\ln x)^{2A}.
\end{align*}
Thus we finally obtain:
\begin{align*}
& j(M)\,\ll\,\frac{1}{M}\int_{0.5}^{1 - \frac{\Z c}{\Z \sqrt{\ln
x}}} N^{\sigma}\,M(\ln x)^{A}d\sigma\,\ll\,N^{1 - \frac{\Z c}{\Z
\sqrt{\ln x}}}(\ln x)^{2A}\,\ll\,xe^{-\frac{\Z c}{\Z
2\mathstrut}\sqrt{\ln x}}.
\end{align*}

\vspace{0.2cm}

$4\text{b}^{\circ}$. The case of <<large>> $M$: $(\ln x)^{A}<M\le
\sqrt{x}(\ln x)^{-A}$. Let us divide the domain of $\gamma$ into
segments $U\!<\!|\gamma|\le U_{1}$ where
\[
U_{1}\le 2U, \quad U_{1}\,\le\,T(q_{1},\chi_{1})\,\le\,2M(\ln
x)^{C}\,=\,T_{1}.
\]
Thus we have
\begin{align*}
& r_{4}(\sigma,M)\,\ll\,\prsum\limits_{U\le
T_{1}}\frac{1}{U}\,r_{4}(\sigma,M,U),\\
& r_{4}(\sigma, M,U)\,=\,\sum\limits_{M<q_{1}\le
M_{1}}\;\stsum\limits_{\chi_{1} \mmod
q_{1}}\;\sum\limits_{\substack{U<|\gamma|\le U_{1} \\ \beta >
0.5}}g_{\chi_{1}}(\varrho,\sigma)|L(\sigma+i\gamma,\chi_{1})|^{\frac{\Z
1}{\Z 2\mathstrut}}.
\end{align*}
Applying H\"{o}lder's inequality, we find
\begin{align*}
& r_{4}^{4}(\sigma, M,U)\,\ll\,\Bigl(\sum\limits_{M<q_{1}\le
M_{1}}\;\stsum\limits_{\chi_{1} \mmod
q_{1}}\;\sum\limits_{\substack{U<|\gamma|\le U_{1} \\ \beta >
0.5}}g_{\chi_{1}}\Bigr)^{3}\times \\
& \times\,\Bigl(\sum\limits_{M<q_{1}\le
M_{1}}\;\stsum\limits_{\chi_{1} \mmod
q_{1}}\;\sum\limits_{\substack{U<|\gamma|\le U_{1} \\ \beta >
0.5}}|L(\sigma + i\gamma,\chi_{1})|\Bigr)\,=\,r_{5}\cdot r_{6}^{3},
\end{align*}
where the notations $r_{5}, r_{6}$ are obvious.

First we get:
\[
r_{5}\,\ll\,\sum\limits_{M<q_{1}\le M_{1}}\;\stsum\limits_{\chi_{1}
\mmod q_{1}}\;\sum\limits_{|\gamma|\le U_{1}}|L(\sigma +
i\gamma,\chi_{1})|^{2}.
\]
Let us divide the domain of $\gamma$ into the segments
\[
2n\,\le\,\gamma\,<\,2n+1,\quad n = 0,\pm 1,\pm 2,\ldots,
\]
and
\[
2n+1\,\le\,\gamma\,<\,2(n+1),\quad n = 0,\pm 1,\pm 2,\ldots\;.
\]
Thus the sum $r_{5}$ splits into sums $r_{5}'$ and $r_{5}''$.
Further, let us sort  all the ordinates $\gamma$ in each segment
$2n\le \gamma < 2n+1$ in increasing order:
\[
2n\le \gamma^{(1)}\le\gamma^{(2)}\le
\ldots\,\le\,\gamma^{(s)}\,<\,2n+1.
\]
It's obvious that $s = O(\ln{q_{1}(|n|+1)})$. Now we place the terms
of $r_{5}'$ that corresponds to the first ordinates $\gamma^{(1)}$
into sum $r_{5}^{(1)}$. The terms corresponding to the second
ordinates $\gamma^{(2)}$ are placed into sum $r_{5}^{(2)}$ and so
on. Thus $r_{5}'$ splits into $s_{0} = O(\ln x)$ sums $r_{5}^{(s)}$,
$s = 1,2,\ldots, s_{0}$. The ordinates $\gamma, \gamma'$ that
corresponds to the neighbouring summands in the sum $r_{5}^{{s}}$
satisfy the condition $|\gamma - \gamma'|>1$. Finally, we apply the
same transformation to the sum $r_{5}''$. Thus the sum $r_{5}$
splits into $\le 2s_{0} = O(\ln x)$ sums $r$ of the following type:
\[
r\,=\,\sum\limits_{M<q_{1}\le M_{1}}\;\stsum\limits_{\chi_{1} \mmod
q_{1}}\prsum\limits_{|\gamma|\le U}|L(\sigma+i\gamma,\chi_{1})|^{2}
\]
where the prime sign means the summation over the <<rarefied>>
ordinates $\gamma$.

Taking $Y = M_{1}(T_{1}+1)$ we have $Y\ge q_{1}(|\gamma|+1)/\pi$ for
all $q_{1}$ and $\gamma$ under considering. Then Lemma 3 implies
that
\begin{align*}
& r\,\ll\, \sum\limits_{M<q_{1}\le M_{1}}\;\stsum\limits_{\chi_{1}
\mmod q_{1}}\prsum\limits_{|\gamma|\le
U}\Bigl(\Bigl|\sum\limits_{n\le Y}\frac{\chi_{1}(n)}{n^{\sigma +
i\gamma}}\Bigr|^{2}\,+\,1\Bigr)\,\ll \\
& \ll\,M^{2}U\,+\,\sum\limits_{M<q_{1}\le
M_{1}}\;\stsum\limits_{\chi_{1} \mmod
q_{1}}\;\prsum\limits_{|\gamma|\le U}\biggl|\sum\limits_{n\le
Y}\frac{\chi_{1}(n)}{n^{\sigma + i\gamma}}\biggr|^{2}.
\end{align*}
The application of Lemma 1 yields:
\[
\prsum\limits_{|\gamma|\le U}\biggl|\sum\limits_{n\le
Y}\frac{\chi_{1}(n)}{n^{\sigma +
i\gamma}}\biggr|^{2}\,\le\,j_{1}\,+\,2\sqrt{j_{1}\,j_{2}},
\]
where
\begin{align*}
& j_{1}\,=\,\int_{-(U_{1}+1)}^{U_{1}+1}\biggl|\sum\limits_{n\le
Y}\frac{\chi_{1}(n)}{n^{\sigma + it}}\biggr|^{2}dt, \quad
j_{2}\,=\,\int_{-(U_{1}+1)}^{U_{1}+1}\biggl|\sum\limits_{n\le
Y}\frac{\chi_{1}(n)\ln n}{n^{\sigma + it}}\biggr|^{2}dt.
\end{align*}
Since $U_{1}+1\le 2U+1\le 3U$, $2\sqrt{j_{1}j_{2}}\le j_{1}+j_{2}$,
by Lemma 2 we obtain
\begin{align*}
& r \,\ll\,M^{2}U\,+\,\sum\limits_{M<q_{1}\le
M_{1}}\;\stsum\limits_{\chi_{1} \mmod q_{1}}(j_{1}\,+\,j_{2})\,\ll\\
& \ll\, M^{2}U\,+\,\int_{-3U}^{3U}\sum\limits_{M<q_{1}\le
M_{1}}\;\stsum\limits_{\chi_{1} \mmod
q_{1}}\biggl(\biggl|\sum\limits_{n\le Y}\frac{\chi_{1}(n)}{n^{\sigma
+ it}}\biggr|^{2}\,+\,\biggl|\sum\limits_{n\le
Y}\frac{\chi_{1}(n)\ln n}{n^{\sigma +
it}}\biggr|^{2}\biggr)dt\,\ll\\
&
\ll\,M^{2}U\,+\,\int_{-3U}^{3U}\bigl(M^{2}+Y\bigr)\Bigl(\sum\limits_{n\le
Y}\frac{1}{n^{2\sigma}}\,+\,\sum\limits_{n\le
Y}\frac{\ln^{2}n}{n^{2\sigma}}\Bigr)dt\,\ll\\
& \ll\, M^{2}U\,+\,\int_{-3U}^{3U}M^{2}(\ln
x)^{C+3}\,dt\,\ll\,M^{2}U(\ln x)^{C+3}
\end{align*}
and therefore
\[
r_{4}(\sigma, M,U)\,\ll\,M^{2}U(\ln x)^{C+4}.
\]
Further, using Lemma 6, we obtain
\begin{align*}
& r_{6}\,\le\,\sum\limits_{q_{1}\le M_{1}}\;\stsum\limits_{\chi_{1}
\mmod q_{1}}\;\sum\limits_{|\gamma|\le
U_{1}}g_{\chi_{1}}(\varrho,\sigma)\,=\,\sum\limits_{q_{1}\le
M_{1}}\;\stsum\limits_{\chi_{1} \mmod q_{1}}N(\sigma; U_{1},
\chi_{1})\;\ll\\
&\ll\, \bigl(M^{2}U\bigr)^{\vth(\sigma)}(\ln x)^{14}.
\end{align*}
Therefore,
\begin{align*}
& r_{4}^{4}(\sigma, M,U)\,\ll\,M^{2}U(\ln
x)^{C+4}\bigl(M^{2}U\bigr)^{3\vth(\sigma)}(\ln
x)^{42}\,\ll\,\bigl(M^{2}U\bigr)^{1+3\vth(\sigma)}(\ln x)^{C+46},\\
& r_{4}(\sigma, M,U)\,\ll\,\bigl(M^{2}U\bigr)^{\frac{\Z 1}{\Z
4\mathstrut}\,+\,\frac{\Z 3}{\Z 4\mathstrut}\vth(\sigma)}(\ln
x)^{\frac{\Z 1}{\Z 4}C +\frac{\Z 23}{\Z 2\mathstrut}},\\
& r_{4}(\sigma,M)\,\ll\,\prsum\limits_{U\le T_{1}}U^{-1}\cdot
\bigl(M^{2}U\bigr)^{\frac{\Z 1}{\Z 4\mathstrut}\,+\,\frac{\Z 3}{\Z
4\mathstrut}\vth(\sigma)}(\ln x)^{\frac{\Z 1}{\Z 4}C +\frac{\Z
23}{\Z 2\mathstrut}}\,\ll \\
& \ll\,(\ln x)^{\frac{\Z 1}{\Z 4}C +\frac{\Z 23}{\Z
2\mathstrut}}\prsum\limits_{U\le T_{1}}U^{\frac{\Z 3}{\Z
4\mathstrut}(\vth(\sigma) - 1)}M^{\frac{\Z 1}{\Z 2\mathstrut} +
\frac{\Z 3}{\Z 2\mathstrut}\vth(\sigma)}.
\end{align*}
Since $\vth(\sigma)\le 1$ for $\tfrac{1}{2}\le \sigma \le 1$, it
follows that
\[
r_{4}(\sigma, M)\,\ll\,M^{\frac{\Z 1}{\Z 2\mathstrut} + \frac{\Z
3}{\Z 2\mathstrut}\vth(\sigma)}(\ln x)^{\frac{\Z 1}{\Z 4}C +\frac{\Z
25}{\Z 2\mathstrut}}
\]
and
\[
j(M)\,=\,M^{-1}\int_{\sigma_{1}}^{1}N^{\sigma}r_{4}(\sigma,M)d\sigma\,\ll\,(\ln
x)^{\frac{\Z 1}{\Z 4}C +\frac{\Z 25}{\Z
2\mathstrut}}\int_{0.5}^{1}\psi(\sigma)\,d\sigma,
\]
where
\[
\psi(\sigma)\,=\,x^{\sigma}M^{u(\sigma)},\quad
u(\sigma)\,=\,-\tfrac{1}{2}\,+\,\tfrac{3}{2}\vth(\sigma).
\]

Let us consider several cases.

$\text{i}^{\circ}$. For $\tfrac{1}{2}\le\sigma\le \tfrac{4}{5}$ we
obviously have
\[
u(\sigma)\,=\,\frac{1}{4}\,+\,\frac{3(4-5\sigma)}{4(2-\sigma)}\,\ge\,\frac{1}{4}.
\]
The condition $N\le \sqrt{x}(\ln x)^{-A}$ implies that
\[
\psi(\sigma)\,\le\, x^{\sigma}\bigl(\sqrt{x}(\ln
x)^{-A}\bigr)^{u(\sigma)}\,\le\,x^{\sigma\,+\,\frac{\Z 1}{\Z
2\mathstrut}u(\sigma)}(\ln x)^{-\frac{\Z 1}{\Z 4\mathstrut}A}.
\]
Since
\[
\sigma\,+\,\frac{1}{2}\,u(\sigma)\,=\,1\,-\,\frac{\bigl(\sigma -
\tfrac{1}{2}\bigr)^{2}}{2\,-\,\sigma}\,\le\,1,
\]
we get
\[
\psi(\sigma)\,\le\,x(\ln x)^{-\frac{\Z 1}{\Z 4\mathstrut}A}.
\]

$\text{ii}^{\circ}$. Suppose $\tfrac{4}{5}\le \sigma\le
\tfrac{6}{7}$. Then
\[
u(\sigma)\,=\,\frac{6 - 7\sigma}{2\sigma}\,\ge\, 0
\]
and
\[
\psi(\sigma)\,\le\,x^{\sigma}\bigl(\sqrt{x}(\ln
x)^{-A}\bigr)^{u(\sigma)}\,\le\,x^{\sigma\,+\,\frac{\Z 1}{\Z
2\mathstrut}u(\sigma)}.
\]
From the inequality
\[
\frac{d}{d\sigma}\Bigl(\sigma\,+\,\tfrac{1}{2}\,u(\sigma)\Bigr)\,=\,1\,-\,\frac{3}{2\sigma^{2}}\,<\,0
\]
it follows that the function $\sigma + \tfrac{1}{2}u(\sigma)$ is
monotonically decreasing on the segment $\tfrac{4}{5}\le\sigma\le
\tfrac{6}{7}$ and attains it's maximum at a point $\sigma=
\tfrac{4}{5}$:
\[
\sigma\,+\,\tfrac{1}{2}\,u(\sigma)\,\le\,\tfrac{4}{5}\,+\,\tfrac{1}{2}\,u\bigl(\tfrac{4}{5}\bigr)\,=\,\tfrac{37}{40}.
\]
Thus
\[
\psi(\sigma)\,\le\,x^{\frac{\Z 37}{\Z 40\mathstrut}}\,\le\,x(\ln
x)^{-\frac{\Z 1}{\Z 4\mathstrut}A}.
\]

$\text{iii}^{\circ}$. Suppose $\tfrac{6}{7}\le \sigma\le
\tfrac{12}{13}$. Then $u(\sigma)\le 0$ and therefore
\[
\psi(\sigma)\,\le\,x^{\sigma}\,\le\,x^{\frac{\Z 12}{\Z
13\mathstrut}}\,\le\,x(\ln x)^{-\frac{\Z 1}{\Z 4\mathstrut}A}.
\]

$\text{iv}^{\circ}$. Suppose that $\tfrac{12}{13}\le \sigma\le 1$.
Then
\[
u(\sigma)\,=\,-\frac{1}{4}\,-\,\frac{13\sigma-12}{4\sigma}\,\le\,-\frac{1}{4}.
\]
Since $M>(\ln x)^{A}$, it implies that
\[
\psi(\sigma)\,\le\,x^{\sigma}M^{-\frac{\Z 1}{\Z
4\mathstrut}}\,\le\,x(\ln x)^{-\frac{\Z 1}{\Z 4\mathstrut}A}.
\]
Thus the inequality
\[
\psi(\sigma)\,\le\,x(\ln x)^{-\frac{\Z 1}{\Z 4\mathstrut}A}
\]
holds for any $\sigma$ such that $\tfrac{1}{2}\le\sigma\le 1$.
Finally we get
\begin{align*}
& j(M)\,\ll\, x(\ln x)^{-\frac{\Z 1}{\Z 4\mathstrut}A\,+\,\frac{\Z
1}{\Z 4\mathstrut}C\,+\,\frac{\Z 25}{\Z 2\mathstrut}}, \\
& R_{4}\,\ll\, (\ln x)^{\frac{\Z 3}{\Z 2\mathstrut}}(\ln\ln
x)^{2}\Bigl(\prsum\limits_{M\le (\ln x)^{C}}xe^{-\frac{\Z c}{\Z
2\mathstrut}\sqrt{\ln x}}\,+\,\prsum\limits_{(\ln x)^{C}<M\le
Q}x(\ln x)^{-\frac{\Z 1}{\Z 4\mathstrut}A\,+\,\frac{\Z 1}{\Z
4\mathstrut}C\,+\,\frac{\Z 25}{\Z 2\mathstrut}}\Bigr)\\
& \ll\, x(\ln x)^{-\frac{\Z 1}{\Z 4\mathstrut}A\,+\,\frac{\Z 1}{\Z
4\mathstrut}C\,+\,15}(\ln\ln x)^{2}.
\end{align*}

Summing the upper bounds for $R_{j}$, $1\le j\le 4$, we obtain
\begin{align*}
& R\,\ll\,(\ln x)^{\alpha}\Bigl(N(\ln x)^{-C+3}\,+\,N(\ln
x)^{-\frac{\Z 3}{\Z 4\mathstrut}C + \frac{\Z 11}{\Z
4\mathstrut}}(\ln\ln x)^{2}\,+\\
& +\,Q\sqrt{N}(\ln x)^{\frac{\Z 1}{\Z 4\mathstrut}(C+1)}(\ln\ln
x)^{2}\,+\,x(\ln x)^{-\frac{\Z 1}{\Z 4\mathstrut}A\,+\,\frac{\Z
1}{\Z 4\mathstrut}C\,+\,15}(\ln\ln x)^{2}\Bigr)\,\ll \\
& \ll (\ln x)^{\alpha}(\ln\ln x)^{2}\Bigl(x(\ln x)^{-C+3}\,+\,x(\ln
x)^{-\frac{\Z 3}{\Z 4\mathstrut}C + \frac{\Z 11}{\Z
4\mathstrut}}\,+\,x(\ln x)^{-A+\frac{\Z 1}{\Z
4\mathstrut}(C+1)}\,+\\
& +\,x(\ln x)^{-\frac{\Z 1}{\Z 4\mathstrut}A\,+\,\frac{\Z 1}{\Z
4\mathstrut}C\,+\,15}(\ln\ln x)^{2}\Bigr)\,\ll\\
& \ll\, x(\ln x)^{\alpha}(\ln\ln x)^{2}\Bigl((\ln x)^{-\frac{\Z
3}{\Z 4\mathstrut}C + \frac{\Z 11}{\Z 4\mathstrut}}\,+\,(\ln
x)^{-\frac{\Z 1}{\Z 4\mathstrut}A\,+\,\frac{\Z 1}{\Z
4\mathstrut}C\,+\,15}\Bigr).
\end{align*}
Taking $C = \tfrac{1}{4}(A - 49)$, we get the required inequality:
\[
R\,\ll\,x(\ln x)^{\alpha}(\ln\ln x)^{2}(\ln x)^{-\frac{\Z 3}{\Z
16\mathstrut}A\,+\,\frac{\Z 191}{\Z 16\mathstrut}}\,\ll\,x(\ln
x)^{-\bigl(\frac{\Z 3}{\Z 16\mathstrut}A\,-\,12\bigr)}.
\]
If we put $B = \tfrac{3}{16}A - 12$, then $A = \tfrac{16}{3}B+4$.
This completes the proof of the lemma.

\vspace{0.2cm}

\begin{flushleft}
\textbf{4. Main theorem}
\end{flushleft}

\vspace{0.2cm}

By Lemma 12, it remains to estimate $R_{a}(x)$ and find an
asymptotic formula for $E_{a}(x)$. Setting $Q = Xd^{-1} =
\sqrt{x}d^{-1}(\ln x)^{-A}$, $N = x(d\delta)^{-1}$ in Lemma 13, we
obtain
\[
R_{a,d,\delta}\,\ll\,x(\ln x)^{-\bigl(\frac{\Z 3}{\Z
16\mathstrut}A\,-\,12\bigr)},\quad R_{a}(x)\,\ll\,x(\ln
x)^{-\bigl(\frac{\Z 3}{\Z 16\mathstrut}A\,-\,15\bigr)}.
\]
Taking $A = 85\tfrac{1}{3}$, we get the formula
\[
S_{a}(x)\,=\,C(\ln x)\beta(a)E_{a}(x)\,+\,O(x\ln\ln x),
\]
where
\[
E_{a}(x)\,=\,\sum\limits_{n\le x}\frac{e_{a}(n)}{\tau(n)}.
\]
In order to calculate the sum $E_{a}(x)$, let us prove that
$e_{a}(n)$ is multiplicative function for any fixed $a\ge 1$.
Indeed, suppose $(m,n) = 1$. Then there exists a unique
decomposition $a = a_{1}a_{2}a'$ where the factors $a_{1}, a_{2},
a'$ are defined as follows. All prime divisors of $a_{1}$ and
$a_{2}$ are among the sets of prime divisors of $m$ and $n$
respectively, and $(a', mn) = 1$. Then
\[
(a_{1},a_{2})\,=\,(a_{1},a')\,=\,(a_{2},a')\,=\,(a_{1},n)\,=\,(a_{2},m)\,=\,1.
\]
Further, if $d$ divides $(a, mn)$ then $d$ has the form $d_{1}d_{2}$
where $d_{1}|(a,m)$ and $d_{2}| (a,n)$. Thus,
\begin{align*}
& e_{a}(mn)\,=\,\sum\limits_{d|(a,mn)}\frac{\D
\beta\left(\frac{an}{d^{2}}\right)}{\beta(a)}\,=\,\sum\limits_{d_{1}|(a,m)}\sum\limits_{d_{2}|(a,n)}
\frac{\D \beta\left(\frac{a_{1}a_{2}a'mn}{d_{1}^{2}d_{2}^{2}}\right)}{\beta(a_{1}a_{2}a')}\,=\\
& =\,\sum\limits_{d_{1}|(a,m)}\sum\limits_{d_{2}|(a,n)}\frac{\D
\beta\left(\frac{a_{1}m}{d_{1}^{2}}\cdot\frac{a_{2}n}{d_{2}^{2}}
\right)\beta(a')}{\beta(a_{1}a_{2})\beta(a')}.
\end{align*}
Since the numbers $\frac{\D a_{1}m}{\D d_{1}^{2}}$ and $\frac{\D
a_{2}n}{\D d_{2}^{2}}$ are integral and coprime, we have
\begin{align*}
& e_{a}(mn)\,=\,\sum\limits_{d_{1}|(a_{1},m)}\frac{\D
\beta\left(\frac{a_{1}m}{d_{1}^{2}}\right)}{\beta(a_{1})}\sum\limits_{d_{2}|(a_{2},n)}\frac{\D
\beta\left(\frac{a_{2}n}{d_{2}^{2}}\right)}{\beta(a_{2})}\,=\\
& =\,\sum\limits_{d_{1}|(a_{1},m)}\frac{\D
\beta\left(\frac{a_{1}m}{d_{1}^{2}}\right)\beta(a_{2}a')}{\beta(a_{1})\beta(a_{2}a')}\,\sum\limits_{d_{2}|(a_{2},n)}\frac{\D
\beta\left(\frac{a_{2}n}{d_{2}^{2}}\right)\beta(a_{1}a')}{\beta(a_{2})\beta(a_{1}a')}\,=\\
& =\, \sum\limits_{d_{1}|(a_{1},m)}\frac{\D
\beta\left(\frac{am}{d_{1}^{2}}\right)}{\beta(a)}\,\sum\limits_{d_{2}|(a_{2},n)}\frac{\D
\beta\left(\frac{an}{d_{2}^{2}}\right)}{\beta(a)}\,=\,
\sum\limits_{d_{1}|(a,m)}\frac{\D
\beta\left(\frac{am}{d_{1}^{2}}\right)}{\beta(a)}\,\sum\limits_{d_{2}|(a,n)}\frac{\D
\beta\left(\frac{an}{d_{2}^{2}}\right)}{\beta(a)}\,=\\
& =\,e_{a}(m)e_{a}(n).
\end{align*}
Suppose $\RRe s > 1$. By definition, put
\[
F_{a}(s)\,=\,\sum\limits_{n =
1}^{+\infty}\frac{e_{a}(n)}{\tau(n)}\,n^{-s}.
\]
Then
\[
F_{a}(s)\,=\,\prod\limits_{p}F_{a,p}(s),\quad
F_{a,p}(s)\,=\,1\,+\,\sum\limits_{k =
1}^{+\infty}\frac{e_{a}(p^{k})}{k+1}\,p^{-ks}.
\]
Since
\[
e_{a}(n)\,=\,\frac{\beta(an)}{\beta(a)}\,=\,\beta(n)
\]
for the case $(a,n) = 1$, we obtain
\begin{align*}
& F_{a}(s)\,=\,\prod\limits_{p\nmid a}\Bigl(1\,+\,\sum\limits_{k =
1}^{+\infty}\frac{\beta(p^{k})}{k+1}\,p^{-ks}\Bigr)\prod\limits_{p|a}F_{a,p}(s)\,=\\
& =\,\prod\limits_{p\nmid a}\Bigl(1\,+\,\beta(p)\sum\limits_{k =
1}^{+\infty}\frac{p^{-ks}}{k+1}\,\Bigr)\prod\limits_{p|a}F_{a,p}(s)\,=\\
& =\,\prod\limits_{p\nmid a}\Bigl(1-\beta(p) -
\beta(p)p^{s}\ln(1-p^{-s})\Bigr)\prod\limits_{p|a}F_{a,p}(s)\,=\,F(s)\psi_{a}(s),
\end{align*}
where
\begin{align*}
& F(s)\,=\,\prod\limits_{p}\Bigl(1-\beta(p) -
\beta(p)p^{s}\ln(1-p^{-s})\Bigr),\\
& \psi_{a}(s)\,=\,\prod\limits_{p|a}\frac{\D 1\,+\,\sum\limits_{k =
1}^{+\infty}\frac{e_{a}(p^{k})}{k+1}\,p^{-ks}}{\D 1-\beta(p) -
\beta(p)p^{s}\ln(1-p^{-s})}.
\end{align*}
Further, $F(s) = \sqrt{\zeta(s)}\,\Phi(s)$, where
\[
\Phi(s)\,=\,\prod\limits_{p}\Phi_{p}(s),\quad
\Phi_{p}(s)\,=\,\Bigl(1-\frac{1}{p^{s}}\Bigr)^{\frac{\Z 1}{\Z
2\mathstrut}}\bigl(1-\beta(p) - \beta(p)p^{s}\ln(1-p^{-s})\bigr).
\]
Setting
\[
1-\beta(p) -
\beta(p)p^{s}\ln(1-p^{-s})\,=\,1\,+\,\frac{1}{2p^{s}}\,+\,u(s),\quad
\Bigl(1-\frac{1}{p^{s}}\Bigr)^{\frac{\Z 1}{\Z
2\mathstrut}}\,=\,1\,-\,\frac{1}{2p^{s}}\,+\,v(s)
\]
and using the decomposition
\begin{align*}
& \beta(p)\,=\,1\,-\,\frac{p}{p^{2}-p+1}\,=\,1 -
\frac{p(p+1)}{p^{3}+1}\,=\,1\,-\,\frac{1}{p}\Bigl(1\,+\,\frac{1}{p}\Bigr)\Bigl(1\,+\,\frac{1}{p^{3}}\Bigr)^{-1}\,=\\
& =\,1\,+\,\sum\limits_{k =
0}^{+\infty}(-1)^{k+1}\Bigl(\frac{1}{p^{3k+1}}\,+\,\frac{1}{p^{3k+2}}\Bigr),
\end{align*}
we find
\begin{align*}
&
u(s)\,=\,\frac{1}{3p^{2s}}\,+\,\frac{1}{4p^{3s}}\,+\,\frac{1}{5p^{4s}}\,+\,\ldots\,-\,\\
&-\,\biggl(\frac{1}{p}\,+\,\frac{1}{p^{2}}\,-\,\frac{1}{p^{4}}\,-\,\frac{1}{p^{5}}\,+\ldots\biggr)
\biggl(\frac{1}{2p^{s}}\,+\,\frac{1}{3p^{2s}}\,+\,\frac{1}{4p^{3s}}\,+\ldots\biggr),\\
& v(s)\,=\,\sum\limits_{k =
2}^{+\infty}2^{-2k}\binom{2k}{k}\frac{p^{-ks}}{2k-1}.
\end{align*}
Suppose now $\sigma > \tfrac{1}{2}$. Then
\begin{align*}
&
|u(s)|\,\le\\
&\le\,\frac{1}{3p^{2\sigma}}\,+\,\frac{1}{4p^{3\sigma}}\,+\,\frac{1}{5p^{4\sigma}}\,+\,\ldots\,+\,
\biggl(\frac{1}{p}\,+\,\frac{1}{p^{2}}\biggr)\biggl(\frac{1}{2p^{\sigma}}\,+\,\frac{1}{3p^{2\sigma}}\,+\,\frac{1}{4p^{3\sigma}}\,+\ldots\biggr)\,\le\\
&
\le\,\frac{1}{3p^{2\sigma}}\biggl(1\,+\,\frac{1}{p^{\sigma}}\,+\,\frac{1}{p^{2\sigma}}\,+\ldots\biggr)\,+\,\frac{3}{2p}\cdot\frac{1}{2p^{\sigma}}
\biggl(1\,+\,\frac{1}{p^{\sigma}}\,+\,\frac{1}{p^{2\sigma}}\,+\ldots\biggr)\,\le\\
&
\le\,\biggl(\frac{1}{3p^{2\sigma}}\,+\,\frac{3}{4p^{\sigma+1}}\biggr)\,\frac{1}{1-p^{-\sigma}}\,\le\,\frac{\sqrt{2}}{\sqrt{2}-1}
\biggl(\frac{1}{3p^{2\sigma}}\,+\,\frac{3}{4p^{\sigma+1}}\biggr),
\end{align*}
and therefore
\begin{equation*}
|u(s)|\,\le\,\frac{\sqrt{2}}{\sqrt{2}-1}\biggl(\frac{1}{3}\,+\,\frac{3}{4}\biggr)\max{\biggl(\frac{1}{p^{2\sigma}},\frac{1}{p^{\sigma+1}}\biggr)}\,\le
\begin{cases}
3.7p^{-2\sigma}, & \text{if}\;\; \tfrac{1}{2}<\sigma \le 1,\\
3.7p^{-(\sigma + 1)}, & \text{if}\;\; \sigma\ge 1.
\end{cases}
\end{equation*}
Further,
\begin{align*}
& |v(s)|\,\le\,\sum\limits_{k =
2}^{+\infty}2^{-2k}\binom{2k}{k}\frac{p^{-k\sigma}}{2k-1}\,=\,1\,-\,\frac{1}{2p^{\sigma}}\,-\,\sqrt{1\,-\,\frac{1}{p^{\sigma}}}\,=\\
& =\, \frac{\D
\Bigl(1-\frac{1}{2p^{\sigma}}\,-\,\sqrt{1\,-\,p^{-\sigma}}\Bigr)\Bigl(1-\frac{1}{2p^{\sigma}}\,+\,\sqrt{1\,-\,p^{-\sigma}}\Bigr)}{\D
\Bigl(1-\frac{1}{2p^{\sigma}}\,+\,\sqrt{1\,-\,p^{-\sigma}}\Bigr)}\,=\\
& =\,\frac{\D
\Bigl(1\,-\,\frac{1}{2p^{\sigma}}\Bigr)^{2}\,-\,\Bigl(1\,-\,\frac{1}{p^{\sigma}}\Bigr)}{\D
\Bigl(1-\frac{1}{2p^{\sigma}}\Bigr)\,+\,\sqrt{1-p^{-\sigma}}}\,\le\,\frac{\D\frac{1}{4p^{2\sigma}}
}{\D 1 -
\frac{1}{2\sqrt{2}}\,+\,\sqrt{1-\frac{1}{\sqrt{2}}}}\,<\,\frac{1}{4p^{2\sigma}}
\end{align*}
for any $\sigma > \tfrac{1}{2}$. In particular, in the case
$\sigma\ge 1$ we get
\[
|v(s)|\,<\,\frac{1}{4p^{2\sigma}}\,\le\,\frac{1}{4p^{\sigma + 1}}.
\]
Thus we obtain for $\tfrac{1}{2}\le\sigma\le 1$
\begin{align*}
&
|\Phi_{p}(s)|\,=\,\Bigl|1\,-\,\frac{1}{4p^{2s}}\,+\,u(s)\,-\,v(s)\,-\,\frac{u(s)}{2p^{s}}\,-\,\frac{v(s)}{2p^{s}}\,-\,u(s)v(s)\Bigr|\,\le\\
\end{align*}
\begin{align*}
&
\le\,1\,+\,\frac{1}{4p^{2\sigma}}\,+\,\frac{3.7}{p^{2\sigma}}\,+\,\frac{1}{4p^{2\sigma}}\,+\,\frac{1}{2p^{\sigma}}
\Bigl(\frac{3.7}{p^{2\sigma}}\,+\,\frac{1}{4p^{2\sigma}}\Bigr)\,+\,\frac{3.7}{4p^{4\sigma}}\,<\\
&
<\,1\,+\,\frac{1}{2p^{2\sigma}}\,+\,\frac{3.7}{p^{2\sigma}}\,+\,\frac{3.95}{2\sqrt{2}p^{2\sigma}}\,+\,\frac{3.7}{8p^{2\sigma}}\,<\,1\,+\,\frac{7}{p^{2\sigma}}\,<\,
\Bigl(1\,-\,\frac{1}{p^{2\sigma}}\Bigr)^{7}
\end{align*}
and therefore
\[
|\Phi(s)|\,<\,\prod\limits_{p}\Bigl(1\,-\,\frac{1}{p^{2\sigma}}\Bigr)^{7}\,=\,\zeta^{7}(2\sigma)\,<\,\Bigl(\frac{2\sigma}{2\sigma
- 1}\Bigr)^{\!7}\,=\,\bigl(\sigma\,-\,\tfrac{1}{2}\bigr)^{-7}.
\]
Similarly, in the case $\sigma \ge 1$ we have
\[
|\Phi(s)|\,<\,\zeta^{7}(\sigma + 1)\,\le\,\zeta^{7}(2)\,=\,O(1).
\]
Thus we get for $\RRe s > 1$:
\[
F_{a}(s)\,=\,\sqrt{\zeta(s)}\,\Phi_{a}(s),\quad
\Phi_{a}(s)\,=\,F(s)\psi_{a}(s)
\]
where the function $\Phi_{a}(s)$ is regular in the half-plane $\RRe
s > \tfrac{1}{2}$ and obeys the inequality
\[
|\Phi_{a}(s)|\,\ll\,\max{\Bigl(\bigl(\sigma -
\tfrac{1}{2}\bigr)^{-7},\, 1\Bigr)}.
\]
Now Lemma 8 implies that
\[
E_{a}(x)\,=\,\frac{x}{\sqrt{\ln
x}}\Bigl(\frac{\Phi_{a}(1)}{\sqrt{\pi}}\,+\,O\Bigl(\frac{1}{\ln
x}\Bigr)\Bigr)
\]
and therefore
\[
S_{a}(x)\,=\,K(a)x\sqrt{\ln x}\,+\,O(x\ln\ln x)
\]
where
\[
K(a)\,=\,\frac{1}{\sqrt{\pi}}\,C\beta(a)\Phi_{a}(1).
\]
Finally, considering the constant $K(a)$, we get $K(a) = K\cdot
\kappa(a)$, where
\begin{align*}
&
K\,=\,\frac{1}{\sqrt{\pi}}\,C\Phi(1)\,=\,\frac{C}{\sqrt{\pi}}\prod\limits_{p}\sqrt{1-\frac{1}{p}}\;\bigl(1\,-\,\beta(p)\,-\,\beta(p)p\ln(1\,-\,p^{-1})\bigr)\,=\\
&
=\,\frac{1}{\sqrt{\pi}}\prod\limits_{p}\Bigl(1\,+\,\frac{1}{p(p-1)}\Bigr)\sqrt{1-\frac{1}{p}}\;\frac{p\bigl(1\,+\,(p-1)^{2}\ln\frac{\D p}{\D p-1}\bigr)}{p^{2}-p+1}\,=\\
& =\,
\frac{1}{\sqrt{\pi}}\prod\limits_{p}\sqrt{1-\frac{1}{p}}\;\frac{p^{2}-p+1}{p(p-1)}\,\frac{p\bigl(1\,+\,(p-1)^{2}\ln\frac{\D
p}{\D p-1}\bigr)}{p^{2}-p+1}\,=\\
&
=\,\frac{1}{\sqrt{\pi}}\prod\limits_{p}\Bigl(\frac{1}{\sqrt{p(p-1)}}\,+\,\sqrt{1-\frac{1}{p}}\;(p-1)\ln\frac{\D
p }{\D p-1}\Bigr),\\
& \kappa(a)\,=\,\beta(a)\psi_{a}(1)\,=\,
\beta(a)\prod\limits_{p|a}\frac{\D 1 + \sum\limits_{k =
1}^{+\infty}\frac{\D e_{a}(p^{k})}{\D k+1}\,p^{-k}}{\D 1 +
\beta(p)\sum\limits_{k = 1}^{+\infty}\frac{p^{-k}}{k+1}}.
\end{align*}
Theorem is completely proved.

\vspace{0.2cm}

\textbf{Remark.} The above arguments may be applied to the sums of
the following type:
\[
C_{a}\,=\,\sum\limits_{a<n\le x}\frac{\tau(n)}{\tau(n-a)},\quad
S_{a,k}(x)\,=\,\sum\limits_{n\le x}\frac{\tau(n)}{\tau_{k}(n+a)},
\]
where $a\ge 1$ and $\tau_{k}(n)$ denotes the number of solutions of
the equation
\[
x_{1}\ldots x_{k}\,=\,n
\]
in the natural numbers $x_{1},\ldots, x_{k}$. One can obtain the
formula
\[
S_{a,k}(x)\,=\,K(a,k)\,x(\ln x)^{\frac{\Z 1}{\Z
k\mathstrut}}\,+\,O(x\ln\ln x).
\]

\vspace{0.2cm}

\begin{flushleft}
\textbf{5. Calculation of the constants} $\boldsymbol{\kappa(a)}$
\end{flushleft}

\vspace{0.2cm}

First we note that $\kappa(a)$ is the multiplicative function of
$a$. Indeed, suppose $a = a_{1}a_{2}$ where $(a_{1},a_{2}) = 1$, and
let $p$ be a prime divisor of $a_{1}$. Then $(p,a_{2}) = 1$ and for
any $k\ge 1$
\begin{align*}
&
e_{a}(p^{k})\,=\,\frac{1}{\beta(a)}\sum\limits_{d|(a,p^{k})}\beta\Bigl(\frac{ap^{k}}{d^{2\mathstrut}}\Bigr)\,=\,\frac{1}{\beta(a_{1})\beta_{a_{2}}}
\sum\limits_{d|(a_{1},p^{k})}\beta\Bigl(a_{2}\cdot\frac{a_{1}p^{k}}{d^{2}\mathstrut}\Bigr)\,=\\
&
=\,\frac{1}{\beta(a_{1})}\sum\limits_{d|(a_{1},p^{k})}\beta\Bigl(\frac{a_{1}p^{k}}{d^{2}\mathstrut}\Bigr)\,=\,e_{a_{1}}(p^{k}).
\end{align*}
Now it follows that
\begin{align*}
&
\kappa(a_{1}a_{2})\,=\\
&
=\,\beta(a_{1}a_{2})\prod\limits_{p|a_{1}}\frac{1\,+\,\sum\limits_{k
= 1}^{+\infty}\frac{\D e_{a}(p^{k})}{\D k+1}\,p^{-k}}{\D
1-\beta(p)\,-\,\beta(p)p\ln\frac{\D p}{\D
p-1}}\cdot\prod\limits_{p|a_{2}}\frac{1\,+\,\sum\limits_{k =
1}^{+\infty}\frac{\D e_{a}(p^{k})}{\D k+1}\,p^{-k}}{\D
1-\beta(p)\,-\,\beta(p)p\ln\frac{\D p}{\D p-1}}\,=\\
& =\,\beta_{a_{1}}\prod\limits_{p|a_{1}}\frac{1\,+\,\sum\limits_{k =
1}^{+\infty}\frac{\D e_{a_{1}}(p^{k})}{\D k+1}\,p^{-k}}{\D
1-\beta(p)\,-\,\beta(p)p\ln\frac{\D p}{\D p-1}}\cdot
\beta_{a_{2}}\prod\limits_{p|a_{2}}\frac{1\,+\,\sum\limits_{k =
1}^{+\infty}\frac{\D e_{a_{2}}(p^{k})}{\D k+1}\,p^{-k}}{\D
1-\beta(p)\,-\,\beta(p)p\ln\frac{\D p}{\D p-1}}\,=\\
& =\,\kappa(a_{1})\kappa(a_{2}).
\end{align*}
Further, $e_{a}(p^{k}) = \beta(p)$ for $(a,p) = 1$ and
\begin{equation*}
e_{a}(p^{k})\,=\,
\begin{cases}
\min{(k+1,m+1)}, & \text{if}\;\; k\ne m,\\
m\,+\,\frac{\D 1}{\D \beta(p)}, & \text{if}\;\; k = m,
\end{cases}
\end{equation*}
for the case $a = p^{m}a_{1}$, $(a_{1},p) = 1$, $m\ge 1$. Now we
calculate the values of $\kappa(a)$ for the cases $a = p, p^{2},
p^{3}, p^{4}$ ($p$ is prime).

\vspace{0.2cm}

$1^{\circ}$. Suppose $a = p$; then
\[
e_{a}(p)\,=\,1\,+\,\frac{1}{\beta(p)},\quad e_{a}(p^{k})\,=\,2,\quad
k = 2,3,4,\ldots\;.
\]
Therefore,
\begin{align*}
& 1\,+\,\sum\limits_{k =
1}^{+\infty}\frac{e_{a}(p^{k})}{k+1}\,p^{-k}\,=\,1\,+\,\frac{1}{2p}\Bigl(1\,+\,\frac{1}{\beta(p)}\Bigr)\,+\,2\sum\limits_{k
= 2}^{+\infty}\frac{p^{-k}}{k+1}\,=\\
&
=\,2p\ln\frac{p}{p-1}\,-\,1\,-\,\frac{1}{2p}\,+\,\frac{1}{2p\beta(p)}.
\end{align*}
Hence,
\[
\kappa(p)\,=\,\frac{\D 2p\ln\frac{\D p}{\D
p-1}\,-\,1\,-\,\frac{1}{2p}\,+\,\frac{\D 1}{\D 2p\beta(p)}}{\D
p\ln\frac{\D p}{\D p-1}-1+\frac{\D 1}{\D \beta(p)}}.
\]
In particular,
\begin{align*}
& \kappa(2)\,=\,\frac{\D 2\ln 2\,-\,\tfrac{1}{4}}{\D \ln 2  +
1}\,=\,0.671\,113\,754\,\ldots\\
& \kappa(3)\,=\,\frac{\D 2\ln \tfrac{3}{2}\,-\,\tfrac{7}{24}}{\D \ln
\tfrac{3}{2} + \tfrac{1}{4}}\,=\,0.792\,206\,241\,\ldots\\
& \kappa(5)\,=\,\frac{\D 2\ln \tfrac{5}{4}\,-\,\tfrac{31}{160}}{\D
\ln \tfrac{5}{4} + \tfrac{1}{16}}\,=\,0.884\,098\,735\,\ldots\\
& \kappa(7)\,=\,\frac{\D 2\ln \tfrac{7}{6}\,-\,\tfrac{71}{504}}{\D
\ln \tfrac{7}{6} + \tfrac{1}{36}}\,=\,0.920\,297\,714\,\ldots\\
& \kappa(11)\,=\,\frac{\D 2\ln
\tfrac{11}{10}\,-\,\tfrac{199}{2200}}{\D \ln
\tfrac{11}{10} + \tfrac{1}{100}}\,=\,0.951\,150\,347\,\ldots\\
& \kappa(13)\,=\,\frac{\D 2\ln
\tfrac{13}{12}\,-\,\tfrac{287}{3744}}{\D
\ln \tfrac{13}{12} + \tfrac{1}{144}}\,=\,0.959\,100\,63\,\ldots\\
& \kappa(17)\,=\,\frac{\D 2\ln
\tfrac{17}{16}\,-\,\tfrac{511}{8704}}{\D \ln
\tfrac{17}{16} + \tfrac{1}{256}}\,=\,0.969\,157\,895\,\ldots\\
& \kappa(19)\,=\,\frac{\D 2\ln
\tfrac{19}{18}\,-\,\tfrac{647}{12312}}{\D \ln \tfrac{19}{18} +
\tfrac{1}{324}}\,=\,0.972\,537\,955\,\ldots\;.
\end{align*}

\vspace{0.2cm}

$2^{\circ}$. Suppose $a = p^{2}$; then
\begin{align*}
& 1\,+\,\sum\limits_{k =
1}^{+\infty}\frac{e_{a}(p^{k})}{k+1}\,p^{-k}\,=\,3p\ln\frac{p}{p-1}\,-\,2\,-\,\frac{1}{2p}\,-\,\frac{1}{3p^{2}}\,+\,\frac{1}{3p^{2}\beta(p)},\\
& \kappa(p^{2})\,=\,\frac{\D
3p\ln\frac{p}{p-1}\,-\,2\,-\,\frac{1}{2p}\,-\,\frac{1}{3p^{2}}\,+\,\frac{1}{3p^{2}\beta(p)}}{\D
p\ln\frac{\D p}{\D p-1}-1+\frac{\D 1}{\D \beta(p)}}.
\end{align*}
In particular,
\begin{align*}
& \kappa(2^{2})\,=\,\frac{\D 3\ln 2\,-\,\tfrac{25}{24}}{\D \ln 2 +
1}\,=\,0.612\,926\,558\,\ldots \\
& \kappa(3^{2})\,=\,\frac{\D 3\ln
\tfrac{3}{2}\,-\,\tfrac{77}{108}}{\D \ln \tfrac{3}{2} +
\tfrac{1}{4}}\,=\,0.768\,053\,638\,\ldots\;.
\end{align*}

\vspace{0.2cm}

$3^{\circ}$. Suppose $a = p^{3}$; then
\begin{align*}
& 1\,+\,\sum\limits_{k =
1}^{+\infty}\frac{e_{a}(p^{k})}{k+1}\,p^{-k}\,=\,4p\ln\frac{p}{p-1}\,-\,3\,-\,\frac{1}{p}\,-\,\frac{1}{3p^{2}}\,-\,\frac{1}{4p^{3}}\,+\,\frac{1}{4p^{3}\beta(p)},\\
& \kappa(p^{3})\,=\,\frac{\D
4p\ln\frac{p}{p-1}\,-\,3\,-\,\frac{1}{p}\,-\,\frac{1}{3p^{2}}\,-\,\frac{1}{4p^{3}}\,+\,\frac{1}{4p^{3}\beta(p)}}{\D
p\ln\frac{\D p}{\D p-1}-1+\frac{\D 1}{\D \beta(p)}}.
\end{align*}
In particular,
\begin{align*}
& \kappa(2^{3})\,=\,\frac{\D 4\ln 2\,-\,\tfrac{169}{96}}{\D \ln 2 +
1}\,=\,0.597\,805\,121\,\ldots\;.
\end{align*}

\vspace{0.2cm}

$4^{\circ}$. Suppose $a = p^{4}$; then
\begin{align*}
& 1\,+\,\sum\limits_{k =
1}^{+\infty}\frac{e_{a}(p^{k})}{k+1}\,p^{-k}\,=\,5p\ln\frac{p}{p-1}-4\,-\,\frac{3}{2p}\,-\,\frac{2}{3p^{2}}\,-\,\frac{1}{4p^{3}}\,-
\,\frac{1}{5p^{4}}+\,\frac{1}{5p^{4}\beta(p)},\\
& \kappa(p^{3})\,=\,\frac{\D
5p\ln\frac{p}{p-1}\,-\,4\,-\,\frac{3}{2p}\,-\,\frac{2}{3p^{2}}\,-\,\frac{1}{4p^{3}}\,-
\,\frac{1}{5p^{4}}+\,\frac{1}{5p^{4}\beta(p)}}{\D p\ln\frac{\D p}{\D
p-1}-1+\frac{\D 1}{\D \beta(p)}}.
\end{align*}
In particular,
\begin{align*}
& \kappa(2^{4})\,=\,\frac{\D 5\ln 2\,-\,\tfrac{2363}{960}}{\D \ln 2
+ 1}\,=\,0.593\,142\,51\,\ldots\;.
\end{align*}

\end{document}